\documentclass{article}

\usepackage{arxiv}

\usepackage[utf8]{inputenc} 
\usepackage[T1]{fontenc}    
\usepackage{hyperref}       
\usepackage{url}            
\usepackage{booktabs}       
\usepackage{amsfonts}       
\usepackage{nicefrac}       
\usepackage{microtype}      
\usepackage{lipsum}

\usepackage{relsize}
\usepackage{verbatim}
\usepackage{graphicx}
\usepackage{dcolumn}
\usepackage{bm}
\usepackage{float}

\usepackage{graphics}
\usepackage{color}
\usepackage{xcolor}
\usepackage{bm}
\usepackage{subfigure}

\usepackage{soul}

\usepackage{amsmath,amssymb,amsfonts}

\usepackage{algorithm}
\usepackage{algpseudocode}
\algnewcommand\server{\item[\textbf{Server execution:}]}%
\algnewcommand\client{\item[\textbf{ClientUpdate($k,w$):}]}%

\usepackage{enumitem}
\setlist[itemize]{leftmargin=*}
\setlist[enumerate]{leftmargin=*}
\usepackage{amsmath,amssymb,amsfonts}
\usepackage{graphics}
\usepackage{color}
\usepackage{xcolor}
\definecolor{rev}{rgb}{0,0,0}
\definecolor{rev2}{rgb}{0,0,0}

\usepackage{array}
\newcolumntype{P}[1]{>{\centering\arraybackslash}p{#1}}
\usepackage{multirow}
\usepackage{cancel}
\usepackage[capitalise]{cleveref} 
\usepackage{enumitem}

\newcommand*{\tran}{^{\mkern-1.5mu\mathsf{T}}}

\usepackage{cite}

\title{On the dual advantage of placing observations through forward sensitivity analysis}

\author{
  Shady E Ahmed \\
  Advanced Computing, Mathematics and Data Division,\\
  Pacific Northwest National Laboratory,\\
  Richland, WA 99354, USA.\\
  \texttt{shady.ahmed@pnnl.gov}
  \And
  Omer San \\
  School of Mechanical \& Aerospace Engineering,\\
  Oklahoma State University,\\
  Stillwater, OK 74078, USA.\\
  \texttt{osan@okstate.edu} 
  \And
  Sivaramakrishnan Lakshmivarahan \\
  School of Computer Science,\\
  University of Oklahoma,\\
  Norman, OK 73019, USA.\\
  \texttt{varahan@ou.edu}
  \And
  John M Lewis\\
  Desert Research Institute,\\
  Reno, NV 89512, USA.\\
\texttt{john.lewis@dri.edu}
}

\begin{document}
\maketitle

\begin{abstract}
The four-dimensional variational data assimilation methodology for assimilating noisy observations into a deterministic model has been the workhorse of forecasting centers for over three decades. While this method provides a computationally efficient framework for dynamic data assimilation, it is largely silent on the important question concerning the minimum number and placement of observations. To answer this question, we demonstrate the dual advantage of placing the observations where the square of the sensitivity of the model solution with respect to the unknown control variables, called forward sensitivities, attains its maximum. Therefore, we can force the observability Gramian to be of full rank, which in turn guarantees efficient recovery of the optimal values of the control variables, which is the first of the two advantages of this strategy. We further show that the proposed strategy of placing observations has another inherent optimality: the square of the sensitivity of the optimal estimates of the control with respect to the observations (used to obtain these estimates) attains its minimum value, a second advantage that is a direct consequence of the above strategy for placing observations. Our analytical framework and numerical experiments on linear and nonlinear systems confirm the effectiveness of our proposed strategy. 
\end{abstract}


\section{Introduction} \label{sec:intro}
Estimating the unknown parameters and initial/boundary conditions (collectively known as the control) of a dynamic model based on a finite set of noisy observations constitutes an important class of inverse problems of interest in a variety of disciplines \cite{tarantola1987inverse,kalnay2003atmospheric,kaipio2006statistical,lewis2006dynamic,biros2011large,edwards2015regional,asch2016data,vieli2006numerical,houtekamer2016review,bannister2017review,carrassi2018data,geer2021learning,brajard2021combining}. At its core, this class of inverse problems gives \textcolor{rev}{rise} to at least four question levels: (1) observability, (2) algorithmic path to estimate the model control, (3) observation count and placement of observation, and (4) assessing the quality of the estimates. At the first level, the observability question may be stated as follows: given the model dynamics and the observation operator, which relates the state of the dynamics to the observable, under what condition is it possible to estimate the unknowns based on a finite set of noisy observation? In a fundamental paper, Kalman \cite{kalman1960general} first provided a binary yes or no answer to this question based on the rank of the observability matrix derived from the linear model dynamics and the linear observation operator. Since then, this basic rank condition has been extended in several directions, including local and global observability of nonlinear systems \cite{casti1985nonlinear,isidori1985nonlinear}. 
Furthermore, building upon the traditions in stability theory \cite{bellman1953stability}, researchers have proposed the concept of the degree of observability and scalar measures to quantify this degree \cite{brown1966not,kang2009quantitative,kang2012optimal,king2015observability}. 

At the second level, given that a system is observable, the question is: what is an efficient algorithmic pathway to compute the estimates of the unknowns? Within the context of geosciences, a variational framework known as the four-dimensional variational data assimilation (4-D VAR) provides an answer in two steps. In the first step, a cost functional is defined by the weighted sum of squared differences between the actual and model counterpart of the observation at the time of the observation. The gradient of this cost functional (also known as the adjoint gradient) is computed using a forward run of the model and a backward sweep of its adjoint model. This process is referred to as the adjoint method (refer to Chapters 22-25 in Lewis \emph{et al.} \cite{lewis2006dynamic} for details). In the second step, this gradient is used in a minimization algorithm (Chapters 9-12, Lewis \emph{et al.} \cite{lewis2006dynamic}) to find an update of the control. These two steps are repeated until the desired convergence is reached. 

The 4-D VAR based approach has been the workhorse of weather forecasting centers around the world since its introduction \cite{le1986variational}. Despite its grand success in delivering forecast products for public consumption for decades, the 4-D VAR framework works with a given set of observations and is largely silent on the third level of inquiry, namely, what is the minimum number of observations required and their placement in the spatiotemporal domain to maximize the effectiveness of computing the estimates? At this juncture, it is useful to review the impact of the number of observations in general. Thanks to advances in the sensor, wireless communication, mass storage, and powerful computing technologies, we are steadily moving away from data-sparse to data-rich regimes. The ability to sample spatiotemporal fields at very high frequency (resulting from ever decreasing sampling intervals) has resulted in truly large datasets that are several orders of magnitude greater than what was available a decade ago. This growth has resulted in two side effects. First, the data exhibit very high correlations, which implies that more data does not translate into more information. Second, it is computationally demanding to ingest all the available data into an assimilation algorithm. These latter considerations highlight the importance of identifying smaller and independent subsets of observations to be used in estimation, which has resulted in a growing body of literature on ``thinning'' and creating of ``super observations'' \cite{ochotta2005adaptive}.

In Lakshmivarahan \emph{et al.} \cite{lakshmivarahan2020controlling,lakshmivarahan2022observability}, a strategy for answering this third level of question within the 4-D VAR framework was provided for the first time. \textcolor{rev}{An analysis based on the forward sensitivity method (FSM) has shown that the observability Gramian, $G$, is a good approximation to the Hessian of the cost function (it is exact for the linear case)}. In addition, the observability Gramian admits an additive decomposition: $G = \sum_{i=1}^{N}G_{i}$, where $G_{i}$ is the contribution to the overall Gramian at the $i^{th}$ observation time and $N$ is the number of observations. While each of these components is symmetric and positive semi-definite, it is shown in \cite{lakshmivarahan2020controlling,lakshmivarahan2022observability} that we can indeed force the Gramian $G$ to be positive definite by placing the $N = (n+p)$ observations where the squares of the forward sensitivities of the model solution (with respect to $n$ components of the initial condition and $p$ parameters including the boundary conditions) attain their maxima. This strategy avoids flat patches in the cost functional by bounding the norm of the adjoint gradient away from zero. This is the first of the two advantages referred to in the title of this paper. The papers by Lewis \emph{et al.} \cite{lewis2020placement,lewis2022placement} and Ahmed \emph{et al.} \cite{ahmed2020forward, ahmed2022forward} contain several applications of this strategy.

At the fourth question level, the emphasis shifts to analyzing the quality of the resulting estimates of the control variables. There are two primary ways to approach this question. The first is to theoretically quantify the asymptotic distribution of the estimates by letting the number of observations increase without bounds. Within the context of \emph{time series analysis} dealing with linear, discrete time, and stochastic dynamic models of the \emph{auto regressive integrated moving average} (ARIMA) types, there is vast literature relating to the asymptotic analysis of the estimates of the parameters of these models \cite{fuller1976introduction, brockwell2006time}. Likewise, there is a large body of results relating to static nonlinear models \cite{seber1988nonlinear,gallant1987nonlinear}. But to our knowledge, there is no such theory of estimates for large-scale nonlinear dynamic models of interest in the geosciences. In the absence of such a statistical theory, we settle for the next best option---an analytical approach to quantify the sensitivity of the estimates with respect to the observations. This quantity is particularly important when the given dynamic model exhibits high sensitivity to errors in the control. In this paper, we prove that the sensitivity of the estimates with respect to the observations attains their minimum value when we place the observations where the squares of the forward sensitivities are maximized as in \cite{lakshmivarahan2020controlling,lakshmivarahan2022observability}. The positive definite Gramian resulting from the proposed strategy for observation placement determines the control estimates that exhibit the smallest possible sensitivity to observations. This is the second advantage of setting the observation placement strategy using the forward sensitivity analysis.

\subsection{Historical remarks}
The four levels of inquiry described above have natural connections to many different areas of the data assimilation and parameter estimation in dynamical systems literature spanning several decades. First is the fundamental structural identifiability question examined by Bellman and Astrom \cite{bellman1970structural}, where it was shown that certain types of model parametrization do not admit a unique solution to the parameter estimation problem. Furthermore, there is extensive literature on lumped parameter system identification \cite{ljung1999system} and adaptive control \cite{narendra1989stable} dealing with parameter estimation. A recent paper by Villaverde \cite{villaverde2019observability} contains an extensive analysis of observability and identification of nonlinear systems of interest in mathematical biology. 

Within the parlance of distributed parameter systems, especially in the context of process control in chemical engineering, there is an extensive body of results relating to parameter estimation, placement of observations, and analysis of the sensitivity of parameter estimates with respect to observations. We refer to the monograph by Ucinski \cite{ucinski2004optimal} and papers by Ala{\~n}a \cite{alana2010optimal} and Christopher and Fathalla \cite{christopher1999sensitivity} for more details. The results in the present paper share some connections with the developments in Christopher and Fathalla \cite{christopher1999sensitivity} and Ala{\~n}a \cite{alana2010optimal}. Nonetheless, these two papers deal only with parameter estimation for distributed parameter systems, while the present study addresses the mutual estimation of initial conditions and parameters (denoted as control) in the same footing for systems described by ordinary/partial differential equations (ODEs/PDEs). For general treatment of sensitivity-based methods, we refer to the foundational works by Cacuci and co-workers \cite{cacuci2003sensitivity,cacuci2005sensitivity}. Some of the ideas developed in the current study can be also connected to the vast and growing literature on adjoint-based analysis of targeted observations and adaptive observations \cite{baker2000observation,langland2004estimation,daescu2008sensitivity,cardinali2009monitoring,gelaro2009examination,majumdar2016review} as well as ensemble-based approaches to observation impact analysis \cite{ancell2007comparing,liu2008estimating}.

\subsection{Significance}
This paper addresses, for the first time, key questions related to the simultaneous estimate of unknown initial conditions and model parameters in variational data assimilation frameworks. The main questions that we answer correspond to observation count and optimal placement as well as the subsequent quality of the inverse problem solution. We present a mathematically concrete strategy for selecting the observation locations by tracking the model forward sensitivity metrics, namely the forecast sensitivity with respect to the initial condition and to the model parameters. Because these sensitivities can generally be positive or negative, we opt for placing the observations at points where the squares of the forward sensitivities attain their maximum. We support the proposed observation placement strategy with a clear theoretical analysis to understand the consequences of this strategy. \textcolor{rev}{In particular, we demonstrate both theoretically and empirically that the advocated observation placement methodology has two main advantages as follows.}

First, by placing the observations where the forward sensitivities are maximized, we enforce the observability Gramian matrix to be positive definite. We illustrate that this Gramian is closely related to the adjoint gradient of the cost function. Specifically, a positive definite observability Gramian leads to non-zero values of the adjoint gradient norm, which accelerates the optimization algorithm that is used to estimate the unknown initial conditions and model parameters. We highlight that the choice of observation locations controls the shape of the resulting cost function (defined using a weighted sum of the squared differences between the collected measurements and the model counterpart of the observation). Placing the observations using the presented forward sensitivity analysis avoids flat regions of the cost function, which would negatively affect the algorithm convergence.

\textcolor{rev2}{Second, we demonstrate that the optimal estimates of the initial conditions and model parameters are sensitive to the locations where measurements are collected. More specifically, we prove that the sensitivities of those estimates (with respect to the observations themselves) reach their minimum values when the forecast sensitivities (with respect to the control) reach their maximum values.} In other words, by placing the observations where the squares of forecast sensitivities with respect to the unknown control are maximized, we guarantee that the solution of the inverse problem is the most robust to small perturbations in observation values. This second result has significant merit in practice, considering the inevitably noisy sensor data with different levels of uncertainties. We back up the analytical framework with numerical experiments using representative systems exhibiting linear and nonlinear dynamics with unknown initial conditions and model parameters.

\subsection{Organization of the paper}
In Section~\ref{sec:math}, we provide a succinct version of the statement of the inverse problem of interest in this paper. In Section~\ref{sec:gram}, we introduce an optimal observation placement strategy using the forward sensitivity analysis that links the observability Gramian and the adjoint gradient. We lay out the analytical framework for defining the dual advantages of the proposed observation placement strategy in Section~\ref{sec:gram} and Section~\ref{sec:dual}. In the latter, we provide a theoretical analysis for evaluating the sensitivity of the control estimates with respect to the observation. In Section~\ref{sec:res}, we illustrate the theory of observation placement through the forward sensitivity analysis using examples corresponding to scalar linear and nonlinear dynamical systems with two elements of control---the initial condition and model parameterization---\textcolor{rev}{as well as two systems governed by PDEs in one and two spatial dimensions}. We draw concluding remarks about the derived optimality conditions for observation placement and the associated twin advantages in Section~\ref{sec:conc}.

\section{Mathematical Preliminaries} \label{sec:math}
We start with a description of the model. Our primary goal is to bring out the beauty and elegance of the underlying ideas and results using simple, easily understandable models. \textcolor{rev}{Therefore, we first derive the analytical framework and perform the theoretical analysis using a scalar (one-dimensional) dynamical system with a single parameter. After that, we provide extensions and insights regarding the higher-dimensional cases.}

\subsection{Model} \label{sec:model}
Let $\mathcal{I} \subset \mathbb{R}$ and $\mathcal{P} \subset \mathbb{R}$ be two subsets of the real line representing the set of all allowed values of the initial condition $x(0)$ and parameter $\alpha$ of a dynamic model, respectively. The state, $x(t) \in \mathbb{R}$, evolves according to the following nonlinear, time-invariant, ordinary differential equation (ODE):
\begin{equation}
    \dot{x} = f(x,\alpha), \label{eq:model}
\end{equation}
where $f:\mathbb{R}\times\mathbb{R} \to \mathbb{R}$, $x(0) \in \mathcal{I}$, and $\alpha \in \mathcal{P}$. Henceforth, the triplet $(f,\mathcal{I},\mathcal{P})$ denotes a class of models of interest. Let $x(t) = x(t,x_0,\alpha)$ denote the solution of \cref{eq:model}. Because the solution depends on $x(0)$ and $\alpha$, the vector $c = (x_0,\alpha)\tran \in \mathbb{R}^2$ is called the control and in the rest of the paper, we use $x(t)$ and $x(t,c)$ interchangeably. 

For the FSM theoretical analysis that we carry out, in addition to the existence and uniqueness of the solution of \cref{eq:model}, we are also interested in the smoothness of the solution regarding the existence of continuous (mixed) partial derivatives of $x(t)$ with respect to $x_0$ and $\alpha$ of order $K \ge 1$. In this context, we recall that a function $g: \mathbb{R} \to \mathbb{R}$ belongs to the class $C^k$ if $g$ and its first $K$ derivatives are continuous. According to a theorem in Chapter 2, Section 7 in Arnold~\cite{arnold1992ordinary}, the solution $x(t)$ of \cref{eq:model} belongs to $C^K$ for some integer $K \ge 1$ if $f(x,\alpha)$ belongs to the same class $C^K$. We define $\mathcal{S}(f) = \{ x(t,c) | c\in \mathcal{I} \times \mathcal{P} \}$ as the ensemble of all possible solutions of \cref{eq:model} for a given $f(x,\alpha) \in C^K$, where $x(0)$ and $\alpha$ are varied in $\mathcal{I}$ and $\mathcal{P}$. For example, if $\dot{x} = ax$, $\mathcal{I}=\mathbb{R}^{+}$ (the positive real line), and $\mathcal{P} = [-1, 1]$, then $\mathcal{S}(f=ax)$ is the set of all exponentials of the form $x(t) = e^{at}x(0)$ where $a\in[-1,1]$ and $x(0) \in \mathbb{R}^{+}$. 

The forward sensitivities of the solution at time $t>0$ with respect to the initial condition $x_0$ and parameter $\alpha$ can be defined as follows:
\begin{equation}
    u(t) = \dfrac{\partial x(t)}{\partial x_0}, \quad \text{and} \quad  v(t) = \dfrac{\partial x(t)}{\partial \alpha}.
\end{equation}
By differentiating \cref{eq:model} with respect to $x(0)$ and $\alpha$, it can be verified that $u(t)$ and $v(t)$ evolve according to linear non-autonomous systems given by
\begin{equation}
    \begin{aligned}
        \dot{u}(t) &= D_f(x(t)) u(t), \quad u(0) = 1,\\
        \dot{v}(t) &= D_f(x(t)) v(t) +  D_f^{\alpha}(x(t)), \quad v(0) = 0,
    \end{aligned} \label{eq:sens}
\end{equation}
where
\begin{equation}
    D_f(x) = \dfrac{\partial f}{\partial x}, \quad \text{and} \quad D_f^{\alpha}(x) = \dfrac{\partial f}{\partial \alpha}.
\end{equation}
We refer to \cite{lakshmivarahan2010forward,lakshmivarahan2017forecast} for the derivation of the forward sensitivity method. By solving \cref{eq:model,eq:sens}, we can compute the evolution of the solution $x(t)$ and its sensitivities with respect to time.

\subsection{Observations} \label{sec:obs}
It is assumed that there is an underlying physical process and that the model $(f, \mathcal{I}, \mathcal{P})$ defined above is faithful to this process in the sense that all the large features of the process are captured by the model. Let $\bar{x}(t)$ for $t \ge 0$ be the true (yet unknown) state of the system under consideration. \textcolor{rev}{It is often the case that we may not be able to observe $\bar{x}(t)$ but only a certain function of it at discrete points in time and space. In addition, observations get corrupted by device errors and measurement noise that is usually modeled as additive Gaussian noise with zero mean and known variance $\sigma^2$.}

Let $h:\mathbb{R} \to \mathbb{R}$, where $h\in C^K$, be the function that maps the true state $\bar{x}(k)$ into the observables $z(k)$ as follows:
\begin{equation}
    z(k) = h(\bar{x}(k)) + \eta(k), \label{eq:obs}
\end{equation}
where $k$ is the discrete time index and $\eta(k)$ is white noise, meaning that the observation $z(k)$ contains noisy information about the (unknown) true state. Finally, we define $\mathcal{O} = \{z(k_i) | 1 \le i \le N \}$ as the set of $N$ observations at discrete times given by $0 \le k_1 \le k_2 \le \dots \le k_N$.

\subsection{Statement of the problem} \label{sec:prob}
Given the model, \cref{eq:model}, and the set $\mathcal{O}$ of observations, our goal is to find the initial condition $x(0) \in \mathcal{I}$ and parameter $\alpha \in \mathcal{P}$ such that the solution $x(t)$ of \cref{eq:model} starting with $c = (x_0, \alpha)\tran$ minimizes the weighted sum of squared error between $z(k)$ and $h(x(k))$. To this end, we define a cost functional as follows:
\begin{equation}
    J(c) = \dfrac{1}{2\sigma^2} \sum_{i=1}^{N} [z(k_i) - h(x(k_i))]^2.
\end{equation}
The problem of minimizing $J(c)$ can be solved by computing the gradient $\nabla J(c) \in \mathbb{R}^2$, called the adjoint gradient, and the minimizer can be sought by using the adjoint gradient in an iterative minimization algorithm \cite{lewis2006dynamic}.

\section{Observability Gramian and Adjoint Gradient} \label{sec:gram}
We first examine the fine structure of the adjoint gradient and its dependence on the observability Gramian using the forward sensitivity analysis \cite{lakshmivarahan2020controlling,lakshmivarahan2022observability}. Let $x(t)$ and $\bar{x}(t)$ be the solution of \cref{eq:model} starting from $c=(x_0,\alpha)\tran$ and \textcolor{rev}{$\bar{c}=(\bar{x}_0,\bar{\alpha})\tran$}, respectively, where the perturbations in the initial conditions and the parameter are given by:
\begin{equation}
    \delta x(0) = \bar{x}_0 - x_0, \qquad \delta \alpha = \bar{\alpha} - \alpha, \qquad  \delta c = (\delta x_0, \delta \alpha)\tran.
\end{equation}
The induced first variation in the solution of \cref{eq:model} resulting from the initial perturbation $\delta c$ in the control $c$ is defined as follows:
\begin{equation}
    \delta x(k) = \bar{x}(k) - x(k). \label{eq:xvar}
\end{equation}
The forecast error, also known as the innovation, can be thus computed as follows:
\begin{equation}
    e(k) = z(k) - h(x(k)) = h(\bar{x}(k)) - h(x(k)) + \eta(k). \label{eq:ek}
\end{equation}
Because $h \in C^K$ and by substituting \cref{eq:xvar} in \cref{eq:ek} and applying first-order Taylor expression, the deterministic component of $e(k)$ can be written as follows:
\begin{equation}
    e(k) = D_h(x(k)) \delta x(k), \label{eq:ek2}
\end{equation}
where $D_h(x) = \dfrac{\partial h}{\partial x}$. Also, recall from first principles $\delta x(k) = u(k) \delta x(0) + v(k) \delta \alpha = F \delta c$, where
\begin{equation}
    F = [u,v] \in \mathbb{R}^{1\times 2}. \label{eq:F}
\end{equation}
From \cref{eq:ek2,eq:F}, we get:
\begin{equation}
    e(k) = D_h(x(k)) F \delta c. \label{eq:ek3}
\end{equation}

Now, let $\delta J$ be the induced variation in $J(c)$ resulting from the initial perturbation $\delta c$ in $c$. Then, it can be verified that \cite{lewis2006dynamic}:
\begin{equation}
    \delta J = -\dfrac{1}{\sigma^2} \sum_{i=1}^{N} e(k_i) D_h(x(k_i)) \delta x(k) = -\dfrac{1}{\sigma^2} \sum_{i=1}^{N} e(k_i) D_h(x(k_i)) F \delta c. \label{eq:dJ}
\end{equation}

Substituting for $e(k)$ using \cref{eq:ek3}, we get an \textcolor{rev}{alternative} expression for $\delta J$ as follows:
\begin{equation}
    \delta J = - \dfrac{1}{\sigma^2} \delta c\tran G \delta c, \label{eq:dJ2}
\end{equation}
where
\begin{equation}
    G = \sum_{i=1}^{N} G_i \in \mathbb{R}^{2\times 2} \quad \text{and } G_i = F\tran D_h^2(x(k_i)) F \in \mathbb{R}^{2\times 2}. \label{eq:G}
\end{equation}

The $2\times 2$ matrix $G$ is known as the observability Gramian, which is always symmetric and, in general, positive semi-definite. The problem of placing observations is determining the minimum number of observations required and where to place them in \textcolor{rev}{order} to make $G$ positive definite. To this and, we rewrite $G$ as
\begin{equation}
    G = \sum_{i=1}^{N} D^2_h(x(k_i)) \begin{bmatrix} u^2(k_i) & u(k_i) v(k_i) \\ v(k_i) u(k_i) & v^2(k_i) \end{bmatrix}.
\end{equation}

Under the regularity assumption that $D_h(x)$ does not vanish along the solution in the set $\mathcal{S}(f)$ in Section~\ref{sec:model} and that $D^2_h(x(k_i))$ is positive, a good rule is to place the observations where $u^2(k_i)$ and $v^2(k_i)$ attain their maximum value. Clearly, in this case, we would need only $N=2$ observations: one at the maximum of $u^2(k)$ and another at the maximum of $v^2(k)$. We note that by placing the observations where $u^2$ and $v^2$ attain their maximum values, the observability Gramian $G$ becomes positive definite, which in turn bounds the norm of the adjoint gradient away from zero (see \cref{eq:dJ2}). By doing so, flat patches of the cost functional are avoided and the convergence of the optimization algorithm is improved.

\section{Natural Consequence of the Placement Strategy: A Theoretical Analysis} \label{sec:dual}
The strategy of placing observations where the square of the forward sensitivities attains a maximum has another natural optimality property. To examine this inherent optimality, we seek an expression for the adjoint gradient starting from the fact that $\delta J = \langle \nabla J(c) , \delta c \rangle$ and at the extremum $\nabla J(c)=0$, where $\langle \cdot , \cdot \rangle$ denotes the inner product. Combining this with \cref{eq:dJ}, we obtain the following conditions for the optimality:
\begin{align}
    g_1 &= \dfrac{\partial J}{\partial x_0} = -\dfrac{1}{\sigma^2} \sum_{i=1}^{N} (z(k_i) - h(x(k_i))) D_h(x(k_i)) U(k_i) = 0, \label{eq:opt1} \\
    g_2 &= \dfrac{\partial J}{\partial \alpha} = -\dfrac{1}{\sigma^2} \sum_{i=1}^{N} (z(k_i) - h(x(k_i))) D_h(x(k_i)) V(k_i) = 0. \label{eq:opt2}
\end{align}
Because the maps $f$ in \cref{eq:model} and $h$ in \cref{eq:obs} are fixed, it can be verified that $g_i = g_i(z,x_0,\alpha)$ is a function of $z,x_0,\alpha$, where $z = (z(k_1), \dots, z(k_N))\tran \in \mathbb{R}^N$. In the following analysis, we consider two cases corresponding to a single observation with a linear observation operator and multiple observations with a generic nonlinear observation operator. Although the former is a subset of the latter, we begin with the simpler case to define the main quantities that enable us to understand and analyze the second advantage of the given observation placement strategy.

\subsection{Case 1: scalar dynamical systems with unknown initial condition}
In the first case, we consider a single observation (i.e., $N = 1$), known model parameter $\alpha$, and unknown initial condition $x_0$. Moreover, we assume that the model state is directly observable (i.e., $h(x)=x$, and $D_h(x) = 1$). Thus, the optimality condition in \cref{eq:opt1} can be reduced to the following:
\begin{equation}
    -\sigma^2 g_1(z,k) = [z(k)-x(k)]u(k) = 0. \label{eq:opt11}
\end{equation}
Because $f$, $h$, and $\alpha$ are invariant, solving \cref{eq:opt11} yields an optimal value of $x_0 = x_0(z, k)$ as a function of $z(k)$ and $k$. We are interested in quantifying the sensitivity of the optimal estimate of $x_0$ with respect to the observation $z$ at time $k$ as follows:
\begin{equation}
    y(k) = \dfrac{\partial x_0(z,k)}{\partial z}.
\end{equation}
To this end, by differentiating both sides of \cref{eq:opt11} with respect to $z = z(k)$, we obtain
\begin{equation}
    0 = -\sigma^2 \dfrac{\partial g_1}{\partial z} = [1-u(k)\dfrac{\partial x_0}{\partial z}] u(k) + [z(k) - x(k)]\dfrac{\partial u(k)}{\partial x_0} \dfrac{\partial x_0}{\partial z}. \label{eq:opt111}
\end{equation}
Solving \cref{eq:opt111}, we obtain the required expression:
\begin{equation}
    y(k) = \dfrac{\partial x_0}{\partial z} = \dfrac{u(k)}{u^2(k) - e(k) \dfrac{\partial u(k)}{\partial x_0}}. \label{eq:y1}
\end{equation}
A little reflection reveals that when $x(t)$ is the optimal trajectory, $e(k)=z(k)-x(k)$ is small and can be set to zero in \cref{eq:y1} to give the following approximation:
\begin{equation}
    y(k) = \dfrac{\partial x_0}{\partial z} \approx \dfrac{u(k)}{u^2(k)} = \dfrac{sign(u(k))}{|u(k)|}. \label{eq:y11}
\end{equation}
\Cref{eq:y11} implies that $y(k)$ is minimum where $u^2(k)$ attains its maximum value. In other words, placing the single observation at time $k^*$, where $u^2(k)$ attains its maximum value, ensures that the sensitivity $y(k^*)$ takes its minimum value. 

For the complementary case when $\alpha$ is not known but $x_0$ is known, the optimality condition in \cref{eq:opt2} becomes:
\begin{equation}
    0 = \sigma^2 g_2(z,k) = [z(k) - x(k)]v(k). \label{eq:opt22}
\end{equation}
Through similar arguments, it can be verified that the optimal $\alpha$ that satisfies \cref{eq:opt22} is a function of $z$ and $k$. We define the sensitivity of the optimal $\alpha$ with respect to the observation $z(k)$ as follows:
\begin{equation}
    w(k) = \dfrac{\partial \alpha(z,k)}{\partial z}.
\end{equation}
Differentiating \cref{eq:opt22} with respect to $z$ and simplifying using the same reasoning used in obtaining \cref{eq:y11}, we get:
\begin{equation}
    w(k) = \dfrac{v(k)}{v^2(k) - e(k) \dfrac{\partial v(k)}{\partial \alpha}} \approx \textcolor{rev2}{\dfrac{sign(v(k))}{|v(k)|}}.
\end{equation}
Thus, $w(k)$ takes on its minimum value at $k = k^*$ where $v^2(k)$ attains its maximum value, implying that the proposed observation placement strategy (i.e., where $v^2(k)$ attains its maximum value) leads to the minimum value of $w(k)$. 

\subsection{Case 2: scalar dynamical systems with unknown parameter and initial condition}
Next, we extend our theoretical analysis to consider the simultaneous estimation of unknown initial conditions and model parameters in the general case of multiple observations and arbitrary observation operators. Solving \cref{eq:opt1,eq:opt2}, it follows that the optimal estimates $x_0 = x_0(z)$ and $\alpha = \alpha(z)$ are functions of the observations, where $z = (z(k_1), \dots, z(k_N))\tran \in \mathbb{R}^N$. Let $x=x(t,x_0(z),\alpha(z))$ be the optimal solution of \cref{eq:model}, starting from the optimal control $c(z) = (x_0 (z), \alpha(z))\tran$. In addition, we suppose that $u(t) = u(t,x_0(z),\alpha(z))$, $v(t) = v(t,x_0(z),\alpha(z))$, and $D_h(x(t))$ are evaluated along the optimal trajectory $x(t)$.

To bring out the key ideas and simplify the algebra, without loss of generality, we set $N=2$ and denote $k_i$ by $i$. Then, the optimality conditions given by \cref{eq:opt1,eq:opt2} for $i = 1,2$ can be rewritten as follows:
\begin{equation}
    \begin{aligned}
        \sigma^2 g_1 &= H_1 = H_{11} + H_{12} = 0,\\
        \sigma^2 g_2 &= H_2 = H_{21} + H_{22} = 0, 
    \end{aligned} \label{eq:opt3}
\end{equation}
where
\begin{align}
        H_{11} &= [z(1) - h(x(1))]D_h(x(1)) u(1), &\quad  H_{12} &= [z(2) - h(x(2))]D_h(x(2)) u(2), \label{eq:H1j} \\
        H_{21} &= [z(1) - h(x(1))]D_h(x(1)) v(1), &\quad  H_{22} &= [z(2) - h(x(2))]D_h(x(2)) v(2), \label{eq:H2j}
\end{align}
where $z(i)$, $x(i)$, $u(i)$, and $v(i)$ are the values of the respective quantities evaluated at time $t= k_i = i$, $1\le i\le 2$. We define the sensitivities of $x_0(z)$ and $\alpha(z)$ with respect to $z_i$ for $i = 1,2$ as follows:
\begin{equation}
    y_i = \dfrac{\partial x_0(z)}{\partial z_i}, \quad \text{and } w_i = \dfrac{\partial \alpha(z)}{\partial z_i}. \label{eq:ysens}
\end{equation}
Our goal is to relate the sensitivities in \cref{eq:ysens} to the forecast sensitivities $u(k)$ and $v(k)$. As illustrated in Case 1 above, this can be accomplished by computing the derivatives of $H_i=0$ in \cref{eq:opt3} with respect to $z_i$ for $1 \le i \le 2$, which, when simplified, gives the sought after relation. The details of this derivation are given in Appendix~A. From \cref{eq:y1w1}--\cref{eq:y2w2} of this appendix, we get the following:
\begin{equation}
    G \begin{bmatrix}y_i \\ w_i \end{bmatrix} = D_h(i) \begin{bmatrix}u_i \\ v_i \end{bmatrix},
\end{equation}
for $i = 1, 2$, where $G=G_1+G_2$, $F(i) = [u(i),v(i)] \in \mathbb{R}^{1\times 2}$. The observability Gramian decomposition $G_i$ can be evaluated as follows:
\begin{equation}
    G_i = F\tran(i) D_h^2(i) F(i)
\end{equation}
where it is assumed that $D_h^2(i) = D_h^2(x(i)) > 0$ along the trajectory of $x$ in \cref{eq:model}. Consequently, the inverse of $G$ controls the behavior of $y_i$ and $w_i$. Setting $D_h^2(i) = d_i$ for simplicity in notation, it can be verified that
\begin{equation}
    G = \begin{bmatrix} 
        d_1^2 u_1^2 + d_2^2 u_2^2    &   d_1^2 u_1v_1 + d_2^2 u_2v_2 \\
        d_1^2 u_1v_1 + d_2^2 u_2v_2  &   d_1^2 v_1^2 + d_2^2 v_2^2
        \end{bmatrix} \label{eq:gramian2}
\end{equation}
and its determinant is given by
\begin{equation}
    |G| = d_1^2 d_2^2 (u_1v_2 - u_2 v_1)^2. \label{eq:det}
\end{equation}
Because the sensitivities can generally be positive or negative, for definiteness, our strategy is based on the square of the sensitivities. Therefore, placing the first observation $z_1$ at the maximum value of $u_1^2$ and the second observation $z_2$ at the maximum value of $v_2^2$ maximizes $|G|$. This in turn reduces the magnitudes of $y_i$ and $w_i$ as desired, per Cramer’s rule. \Cref{table:yw} provides expressions for $y(i)$ and $w(i)$ for $i = 1,2$.

\textcolor{rev}{\paragraph{Remark 4.1 ---} For an increased number of observations, it is crucial to guarantee that the observation placement strategy does not lead to a singular Gramian matrix. In Appendix~B, we prove that the forecast sensitivities with respect to the initial condition and parameter cannot be linearly dependent and that the resulting Gramian is therefore non-singular.}

\begin{table*}[htbp!]
\caption{Expressions for $y(i)$ and $w(i)$ for $i = 1,2$.} \vspace{5pt}
\centering
\begin{tabular}{P{0.05\textwidth} P{0.25\textwidth} P{0.25\textwidth}   }  
\hline
& $y$ & $w$ \smallskip \\
\hline \smallskip\\
$i=1$ & $\dfrac{v_2(u_1v_2-u_2v_1)}{|G|}$ & $\dfrac{u_2(u_2v_1-u_1v_2)}{|G|}$  \medskip \vspace{10pt} \\
$i=2$ & $\dfrac{v_1(u_2v_1-u_1v_2)}{|G|}$ & $\dfrac{u_1(u_1v_2-u_2v_1)}{|G|}$ \medskip \\
\hline
\end{tabular}
\label{table:yw}
\end{table*}

\subsection{Extensions to high-dimensional systems}
\textcolor{rev}{We consider the case for which the state variable $\mathbf{x}$ is represented by a vector (i.e., $\mathbf{x}\in \mathbb{R}^n$, where $n$ is the number of degrees of freedom of the system) and $\boldsymbol{\alpha} \in \mathbb{R}^{p}$. The discrete-time version of the model in \cref{eq:model} for $\mathbf{x}\in \mathbb{R}^n$ can be written as follows:
\begin{equation}
    \mathbf{x}_{k+1} = \boldsymbol{M}(\mathbf{x}_k,\boldsymbol{\alpha}), \label{eq:model_discrete_vector}
\end{equation}
where $\boldsymbol{M}:\mathbb{R}^{n}  \times \mathbb{R}^p \to \mathbb{R}^{n}$ defines the one-time step mapping defined by applying a temporal integration scheme (e.g., $\boldsymbol{M}(x_k,\boldsymbol{\alpha}) = \mathbf{x}_k + (t_{k+1}-t_k)\boldsymbol{f}(\mathbf{x}_k,\boldsymbol{\alpha})$ for a first-order Euler integrator). \Cref{eq:sens} for the dynamics of $\mathbf{u}(t) \in \mathbb{R}^{n\times n}$ and $\mathbf{v}(t) \in \mathbb{R}^{n\times p}$ can be rewritten as follows:
\begin{equation}
    \mathbf{u}_{k+1} = \mathbf{D}_{\mathbf{M}}(k) \mathbf{u}_k, \qquad \text{and} \qquad  \mathbf{v}_{k+1} = \mathbf{D}_{\mathbf{M}}(k) \mathbf{v}_k + \mathbf{D}_{\mathbf{M}}^{\boldsymbol{\alpha}}(\mathbf{x}_k),
\label{eq:sens_discrete_vector}
\end{equation}
where
\begin{equation}
    \bigg[\mathbf{D}_{\mathbf{M}}(k)\bigg]_{i,j} = \dfrac{\partial M_i}{\partial x_j}\bigg|_{\mathbf{x}=\mathbf{x}_k}, \qquad \text{and} \qquad \bigg[\mathbf{D}_{\mathbf{M}}^{\boldsymbol{\alpha}}(k)\bigg]_{i,j} = \dfrac{\partial M_i}{\partial \alpha_j}\bigg|_{\mathbf{x}=\mathbf{x}_k}.
\end{equation}
Without loss of generality, we focus our analysis on estimating the model's initial condition $\mathbf{x}_0$ from a set of measurements, $\mathbf{z} \in \mathbb{R}^m$, (i.e., the control is $\mathbf{c}=\mathbf{x}_0$).}

\subsubsection{Linear dynamics (vector)}
\textcolor{rev}{We first consider a linear dynamical system:
\begin{equation}
    \mathbf{x}_{k+1} = \mathbf{M} \mathbf{x}_k, \label{eq:model_discrete_linear_vector}
\end{equation}
where $\mathbf{M}\in \mathbb{R}^{n\times n}$ is a square matrix with rank $n$. Furthermore, we suppose that the measurement vector $\mathbf{z}\in \mathbb{R}^m$ is related to the state $\mathbf{x}$ by a linear operator $\mathbf{H}\in \mathbb{R}^{m\times n}$ as follows:
\begin{equation}
    \mathbf{z}_k = \mathbf{H} \mathbf{x}_k + \boldsymbol{\eta}_k,
\end{equation}
where $ \boldsymbol{\eta}_k \in \mathbb{R}^m$ represents an additive measurement noise. The cost function for the inverse problem can be written as:
\begin{equation}
    \begin{aligned}
        J(\mathbf{c}) &= \dfrac{1}{2} \sum_{i=1}^{N} [\mathbf{z}_i-\mathbf{H} \mathbf{x}_i]\tran R^{-1} [\mathbf{z}_i-\mathbf{H} \mathbf{x}_i].
    \end{aligned} \label{eq:Jc1}
\end{equation}
From \cref{eq:model_discrete_linear_vector}, we can relate the prediction $\mathbf{x}_i$ to the initial condition $\mathbf{x}_0=\mathbf{c}$ as $\mathbf{x}_i = \mathbf{M}^i \mathbf{x}_0= \mathbf{M}^i \mathbf{c}$ and $\mathbf{x}_i\tran =\mathbf{c}\tran (\mathbf{M}^i)\tran = \mathbf{c}\tran (\mathbf{M}\tran)^i$. We also note that the forward sensitivity of the model forecast with respect to its initial condition can be defined as $\mathbf{u}_i = \mathbf{M}^i$ because the model Jacobian is $\mathbf{D}_{\mathbf{M}}=\mathbf{M}$ and $\mathbf{u}_0=\mathbf{I}$ (the identity matrix). Therefore, \cref{eq:Jc1} and its gradient can be rewritten as:
\begin{align}
    J(\mathbf{c}) &= \dfrac{1}{2} \sum_{i=1}^{N} \bigg[\mathbf{z}_i\tran \mathbf{R}^{-1} \mathbf{z}_i - 2 \mathbf{z}_i\tran \mathbf{R}^{-1} \mathbf{H} \mathbf{M}^i \mathbf{c} + \mathbf{c}\tran \big((\mathbf{M}\tran)^i \mathbf{H}\tran \mathbf{R}^{-1}\mathbf{H} \mathbf{M}^i \big) \mathbf{c}\bigg], \\
    \nabla_{\mathbf{c}} J &= \sum_{i=1}^{N} [(\mathbf{M}\tran)^i \mathbf{H}\tran \mathbf{R}^{-1} [\mathbf{H} \mathbf{M}^i \mathbf{c} - \mathbf{z}_i]. \label{eq:gradJc1}
\end{align}
Setting the gradient $\nabla_{\mathbf{c}} J$ in \cref{eq:gradJc1} to zero, we get an optimal expression for $\mathbf{c}$ as follows:
\begin{equation}
    \sum_{i=1}^{N} (\mathbf{M}\tran)^i \mathbf{H}\tran \mathbf{R}^{-1} \mathbf{H} \mathbf{M}^i \mathbf{c}  = \sum_{i=1}^{N} (\mathbf{M}\tran)^i \mathbf{H}\tran \mathbf{R}^{-1} \mathbf{z}_i.
\end{equation}
The observability Gramian is now defined as $\mathbf{G} = \sum_{i=1}^{N} \mathbf{G}_i$, where $\mathbf{G}_i = (\mathbf{M}\tran)^i \mathbf{H}\tran \mathbf{R}^{-1} \mathbf{H} \mathbf{M}^i$. Thus, the optimal $\mathbf{c}$ is given by
\begin{equation}
    \mathbf{c} = \mathbf{G}^{-1} \sum_{i=1}^{N} (\mathbf{M}\tran)^i \mathbf{H}\tran \mathbf{R}^{-1} \mathbf{z}_i. \label{eq:control_linear}
\end{equation}
Differentiating both sides of \cref{eq:control_linear} with respect to $\mathbf{z}_i$ yields the following expression for the sensitivity of the optimal estimate of the control $\mathbf{c}$ to the measurements:
\begin{equation}
    \dfrac{\partial \mathbf{c}}{\partial \mathbf{z}_i} = \mathbf{G}^{-1} (\mathbf{M}\tran)^i \mathbf{H}\tran \mathbf{R}^{-1}.
\end{equation}
Maximizing the square of the forward sensitivity of the model forecast with respect to the control yields a Gramian matrix $\mathbf{G}$ with the largest determinant, which minimizes the sensitivity of the optimal estimate of the control to the actual measurements.}

\textcolor{rev}{We note that the minimum number of observation times $N$ to yield a well-posed inverse problem is defined as $N=\dfrac{n}{m}$. However, we highlight that due to the existence of measurement noise and correlations between measurement components, a regularization is often imposed for the solution of the inverse problem. For the special case of $m=n$, we place a single measurement $\mathbf{z}_i \in \mathbb{R}^n$ at arbitrary time $i$. Assuming $\mathbf{H}=\mathbf{I}$ and $\mathbf{R}^{-1}= \dfrac{1}{\sigma^2} \mathbf{I}$, we get the following expressions for the Gramian and the sensitivity of the optimal estimate of the control $\mathbf{c}$ with respect to the measurement $\mathbf{z}_i$:
\begin{align}
    \mathbf{G} = \dfrac{1}{\sigma^2} (\mathbf{M}\tran)^i \mathbf{M}^i \qquad \text{and} \qquad
    \dfrac{\partial \mathbf{c}}{\partial \mathbf{z}_i} = \mathbf{M}^{-i}.\label{eq:dcdz_linear}
\end{align}
As noted above, the forecast sensitivity with respect to the initial conditions can be written as $\mathbf{u}_i = \mathbf{M}^i$ in the linear dynamics case. Therefore, the sensitivity $\dfrac{\partial \mathbf{c}}{\partial \mathbf{z}_i}$ is inversely proportional to the sensitivity of $\mathbf{x}_i$ with respect to $\mathbf{c}$.}

\textcolor{rev}{\paragraph{Remark 4.2 ---} The system matrix $\mathbf{M}$ in \cref{eq:model_discrete_linear_vector} is a real matrix and it may or may not be diagonalizable. For the latter, it is similar to the Jordan normal form. For the case where $\mathbf{M}$ is diagonalizable, there exists a non-singular matrix $P$ such that
\begin{equation}
    \Lambda = \mathbf{P}^{-1} \mathbf{M} \mathbf{P},
\end{equation}
where $\Lambda$ is a diagonal matrix of the real eigenvalues of $\mathbf{M}$. Therefore, $\mathbf{M}^i=\mathbf{P} \Lambda^i \mathbf{P}^{-1}$, where $\Lambda^i = \text{Diag}(\lambda_1^i, \lambda_2^i, \dots, \lambda_n^i)$, where $\lambda_1 \ge \lambda_2 \ge \dots \ge \lambda_n$. In this case, \cref{eq:dcdz_linear} can be rewritten as follows:
\begin{equation}
    \mathbf{P}^{-1} \bigg[\dfrac{\partial \mathbf{c}}{\partial \mathbf{z}_i} \bigg] \mathbf{P} =\lambda^{-i}. \label{eq:dc_diagonal}
\end{equation}
Therefore, if $|\lambda_k|>1$, $\lambda_k^{i}$ is large (as well as the square of forward sensitivities because $\mathbf{u}_i = \mathbf{M}^i$) and $\dfrac{\partial \mathbf{c}}{\partial \mathbf{z}_i}$ decreases. In addition, \cref{eq:dc_diagonal} shows that the basis $\mathbf{P}$ that diagonalizes $\mathbf{M}$ also diagonalizes the matrix $\bigg[\dfrac{\partial \mathbf{c}}{\partial \mathbf{z}_i}\bigg]_{n\times n}$, which implies that $\bigg[\dfrac{\partial \mathbf{c}}{\partial \mathbf{z}_i}\bigg]_{n\times n}$ is similar to a diagonal matrix.}

\subsubsection{Nonlinear dynamics (vector)}
\textcolor{rev}{We consider the nonlinear case of \cref{eq:model_discrete_linear_vector} as follows:
\begin{equation}
    \begin{aligned}
        \mathbf{x}_{k+1} = M(\mathbf{x}_k), \qquad 
        z_k = h(\mathbf{x}_k) + \eta_k.
    \end{aligned}
\end{equation}
The forecast error is defined as $e_k = z_k - h(\mathbf{x}_k)$, and the cost function can be written as follows:
\begin{equation}
    J(\mathbf{c}) = \dfrac{1}{2} \sum_{i=1}^N e_k\tran R^{-1} e_k = \sum_{i=1}^N J_k(\mathbf{c}),
\end{equation} 
where $J_k(\mathbf{c}) = \dfrac{1}{2} e_k\tran R^{-1} e_k$. To simplify our analysis, we assume $R = \sigma^2 I$, which gives the following:
\begin{equation}
    \begin{aligned}
        J_k(\mathbf{c}) = \dfrac{1}{2\sigma^2} e_k\tran e_k = \dfrac{1}{2\sigma^2} \sum_{j=1}^m (z_{i,j} - h_j(x_i))^2,
    \end{aligned}
\end{equation}
where $z_{i,j}$ is the $j^{th}$ component of the measurement vector $z_i$ at time $t_i$. Thus, the cost function $J(c)$ can be written as follows:
\begin{equation}
    J(c) = \dfrac{1}{2\sigma^2} \sum_{i=1}^N \sum_{j=1}^m (z_{i,j} - h_j(\mathbf{x}_i))^2.
\end{equation}
In other words, $J(\mathbf{c})$ is additive in $J_{k,i} := \dfrac{1}{2\sigma^2} (z_{i,j} - h_j(\mathbf{x}_i))^2 $. Therefore, it is sufficient to investigate $J_{k,i}$ for the analysis of the dual advantages of the proposed forward sensitivity observation placement strategy. The analysis of $J_{k,i}$ is similar to the analysis of $J(c)$ for the scalar case presented in Section~\ref{sec:dual} and Appendix~A.}

\section{Numerical Experiments} \label{sec:res}
In this section, we consolidate the results related to the application of the theoretical analysis carried out in Sections~\ref{sec:gram} and~\ref{sec:dual} to the dynamical system introduced in Section~\ref{sec:math}. \textcolor{rev}{First, we explore the advantages of observation placement using the forward sensitivity method with two test cases for linear and nonlinear scalar dynamical systems.} In both cases, the initial condition and model parameter represent the control variable. We investigate the dependence of the optimal estimate of the sought control variable on the time instants at which measurements are collected. The forward sensitivities (i.e., $u$ and $v$) as well as the sensitivities of the control variable estimates with respect to the different measurements are computed to analyze the effect of observation placement on solving the inverse problem. \textcolor{rev}{After that, we showcase the success of the FSM-based observation placement strategy using two problems governed by PDEs, namely the one-dimensional (1D) Burgers and the two-dimensional (2D) advection diffusion test cases.}

\subsection{Linear dynamics (scalar)}
\begin{figure}[ht!]
    \centering
    \includegraphics[width=0.6\linewidth]{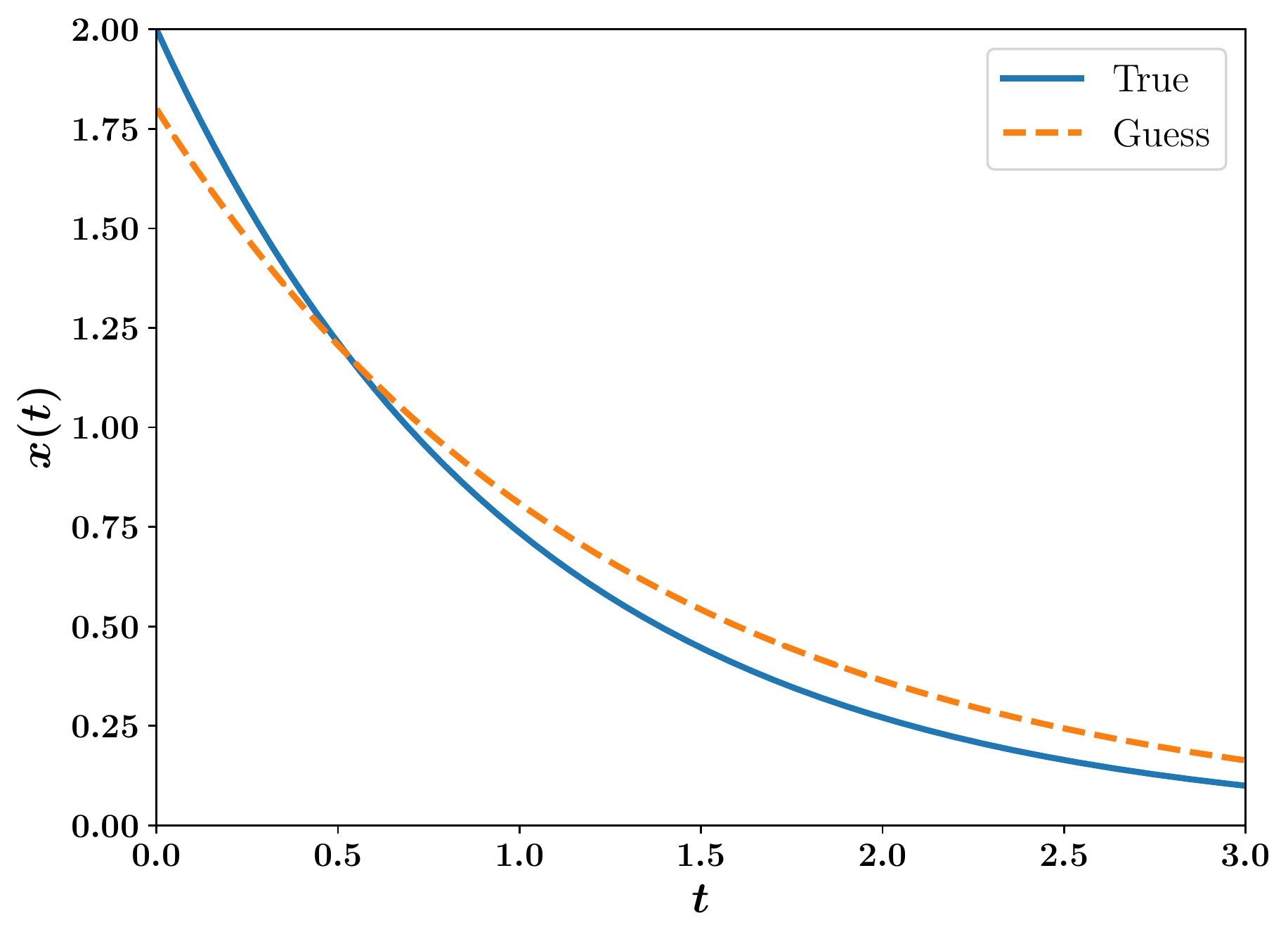}
    \caption{Solution of the decay model using the true and guess controls.}
    \label{fig:xtrue}
\end{figure}

The first example that we consider is an exponential decay model defined as follows:
\begin{equation}
    \dot{x} = ax,
\end{equation}
where $x$ is the model state and $a<0$ is the decay parameter. It can be verified that the analytic solution of this ODE is $x(t) = x_0 e^{at}$, where $x(0)=x_0$ is the initial condition. For our numerical experiments, we set the true controls as $\bar{c}=[\bar{x}_0, \bar{a}]\tran = [2,-1]\tran$. These are assumed to be unknown, and instead the initial guessed values are given by $c=[x_0, a]\tran = [1.8,-0.8]\tran$. The dynamical evolution of the decay in response to the true and guess controls is shown in \cref{fig:xtrue}.

\begin{figure}[ht!]
    \centering
    \includegraphics[width=0.95\linewidth]{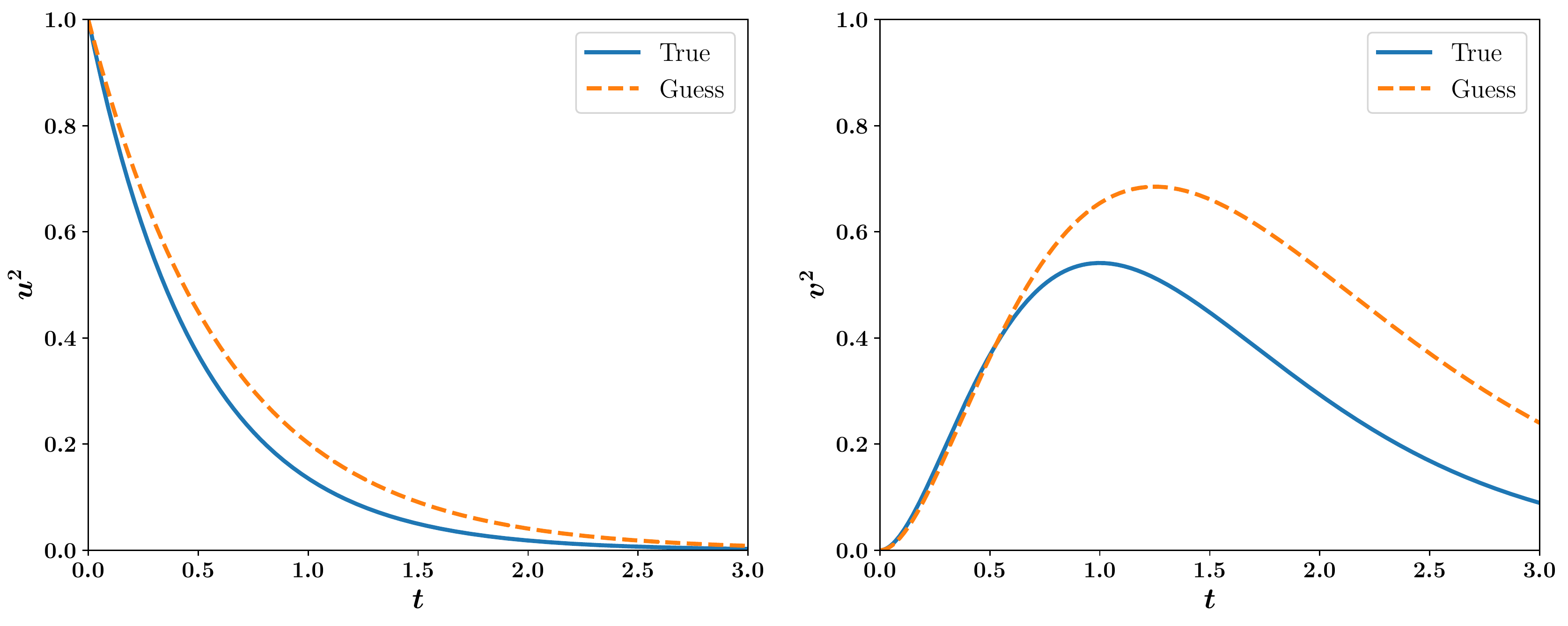}
    \caption{Model sensitivity with respect to (a) initial condition $x_0$ (left) and (b) the decay parameter $a$ (right).}
    \label{fig:sensitivity}
\end{figure}

Synthetic observations are created from the true state with added noise of zero mean and a standard deviation $10\%$ of the model state at times when the observations are made. We follow the proposed methodology for placing the observation where the forward sensitivities attain their maximum values. For the considered exponential decay model, the forward sensitivities of the model predictions with respect to the control variables can be written as follows:
\begin{equation}
    u(t) = e^{at}, \quad v(t) = t e^{at} x_0 = t x(t).
\end{equation}
With $a < 0$, $u(t) = e^{-|a|t}$ attains the maximum at $t=0$ and $v(t) = te^{-|a|t}x$ attains the maximum at $t^*=\dfrac{1}{|a|}$. To avoid a zero value of $v(t)$, we avoid placing observations at $t=0$. Moreover, it can be verified from \cref{eq:det} that $|G| = 0$ if $t_1=t_2$. Therefore, we can place the two observations at $t_1 = \epsilon > 0$ where $\epsilon$ is a small value (as opposed to having the observation coincide with the initial time) and $t_2=t^*=\dfrac{1}{|a|}$. In this case, it can be shown that $G$ is given by
\begin{equation}
    G = e^{-2|a| \epsilon} \begin{bmatrix} 1 & \epsilon x_0 \\ \epsilon x_0 & \epsilon^2 x_0^2 \end{bmatrix} 
      + e^{-2|a| t^*} \begin{bmatrix} 1 & t^* x_0 \\ t^*x_0 & t^{*2} x_0^2 \end{bmatrix},
\end{equation}
which is positive definite. The plots of forecast sensitivities to initial condition $x_0$ and decay parameter $a$ are shown in \cref{fig:sensitivity} for the true and guessed values of control $c$. In order to optimize the solution of the inverse problem, observation sites are chosen where the forecast sensitivities to control reach a maximum (i.e., $u^2(t)$ and $v^2(t)$ exhibit maxima). In particular, we place an observation at $t_1 = k_1 = \epsilon = 0.1$ (to avoid singularity of the Gramian, G, in \cref{eq:det}) and at $t_2 = k_2 = 1.0$.

\begin{figure}[ht!]
    \centering
    \includegraphics[width=0.65\linewidth]{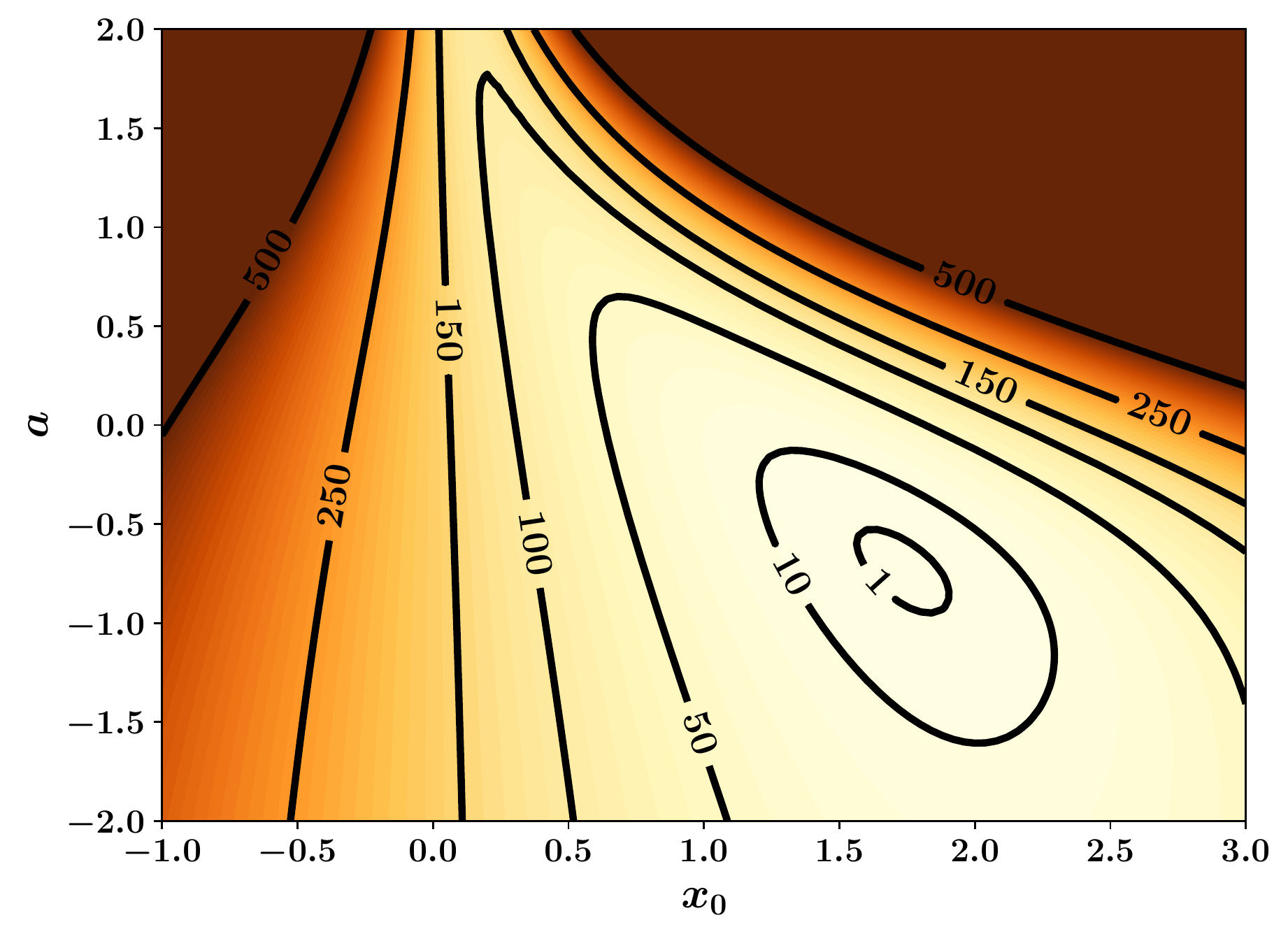}
    \caption{Cost function for exponential decay model.}
    \label{fig:cost}
\end{figure}

The resulting cost function is displayed in \cref{fig:cost} for different control values where the minimum value occurs around the true values of the control. Numerically, the search for the minimum can be made using Newton’s method, leading to optimal estimates of the control values as $c = (x_0, a) = (1.7329,-0.7256)$.

\begin{figure}[ht!]
    \centering
    \includegraphics[width=0.99\linewidth]{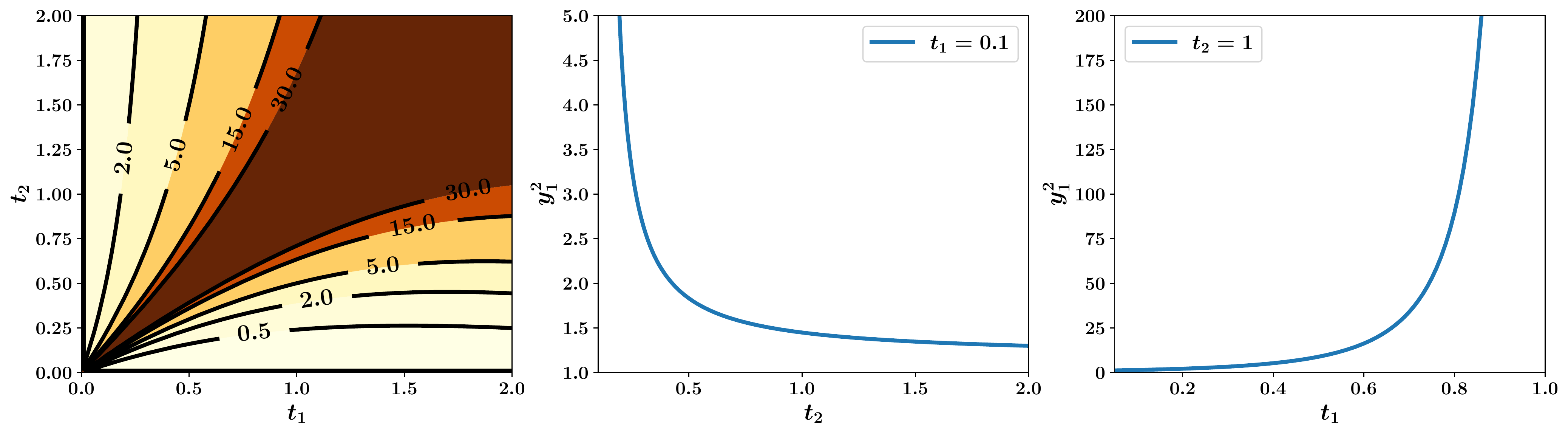}
    \caption{Sensitivity of the optimal estimates of the initial condition $x_0$ to the first observation $z_1$. Left: contour plot of $y_1^2(t)$ at different values of $t_1$ and $t_2$. Middle: $y_1^2(t)$ values at $t_1=0.1$ and different values of $t_2$. Right: $y_1^2(t)$ values at $t_2=1.0$ and different values of $t_1$. Note that minimum values of $y_1^2$ are obtained at $t_1=0.1$ and $t_2=1.0$.}
    \label{fig:y1}
\end{figure}

\begin{figure}[ht!]
    \centering
    \includegraphics[width=0.99\linewidth]{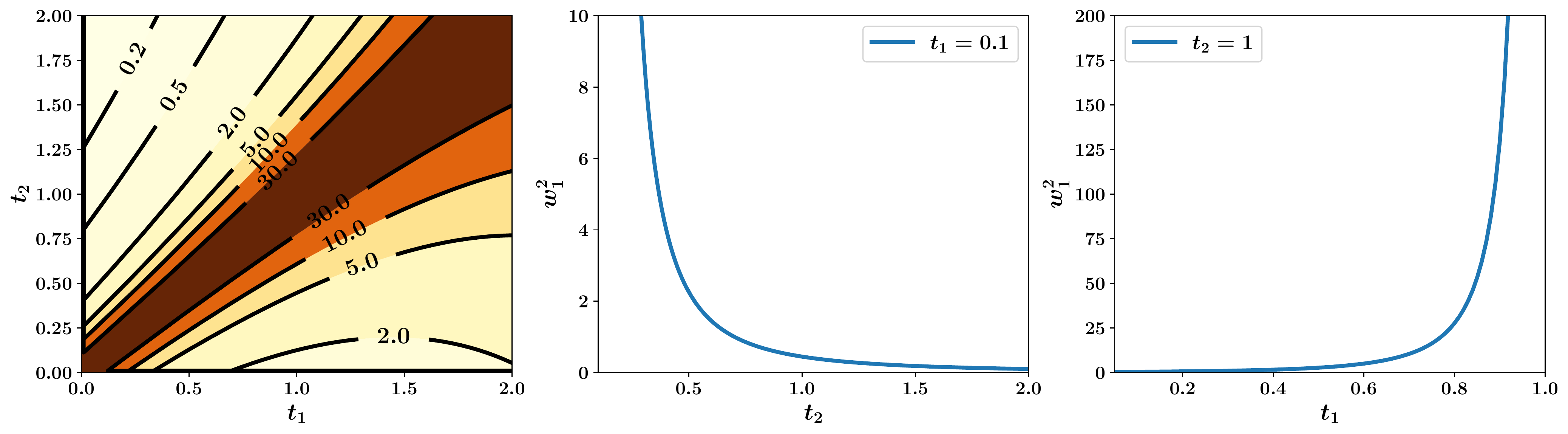}
    \caption{Sensitivity of the optimal estimates of the decay parameter $a$ to the first observation $z_1$. Left: contour plot of $w_1^2(t)$ at different values of $t_1$ and $t_2$. Middle: $y_1^2(t)$ values at $t_1=0.1$ and different values of $t_2$. Right: $w_1^2(t)$ values at $t_2=1.0$ and different values of $t_1$. Note that minimum values of $w_1^2$ are obtained at $t_1=0.1$ and $t_2=1.0$.}
    \label{fig:w1}
\end{figure}

\begin{figure}[ht!]
    \centering
    \includegraphics[width=0.99\linewidth]{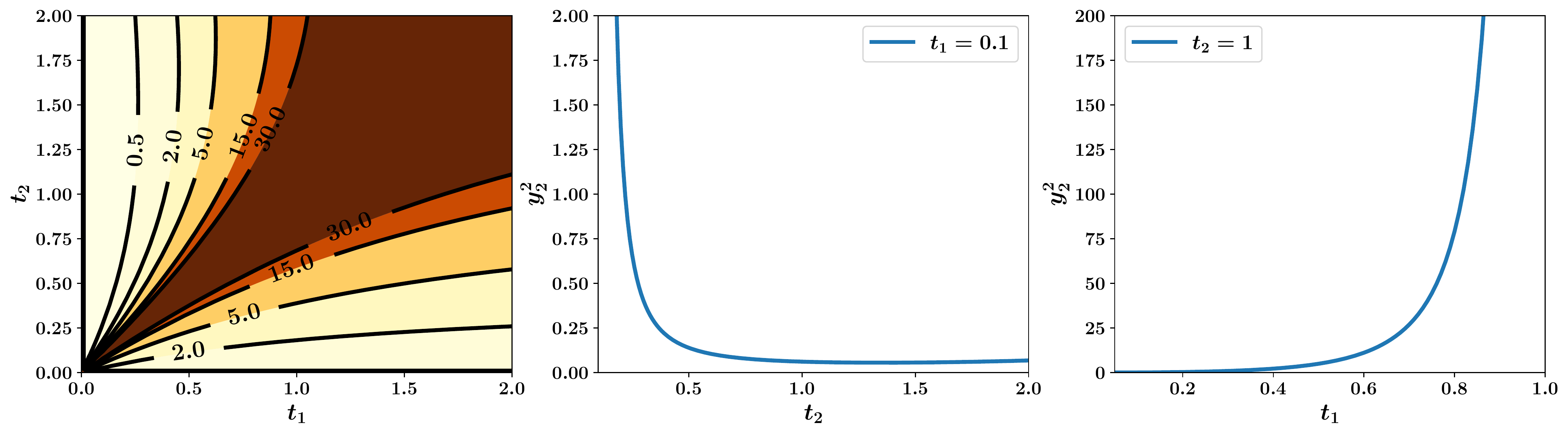}
    \caption{Sensitivity of the optimal estimates of the initial condition $x_0$ to the second observation $z_2$. Left: contour plot of $y_2^2(t)$ at different values of $t_1$ and $t_2$. Middle: $y_2^2(t)$ values at $t_1=0.1$ and different values of $t_2$. Right: $y_1^2(t)$ values at $t_2=1.0$ and different values of $t_1$. Note that minimum values of $y_2^2$ are obtained at $t_1=0.1$ and $t_2=1.0$.}
    \label{fig:y2}
\end{figure}

\begin{figure}[ht!]
    \centering
    \includegraphics[width=0.99\linewidth]{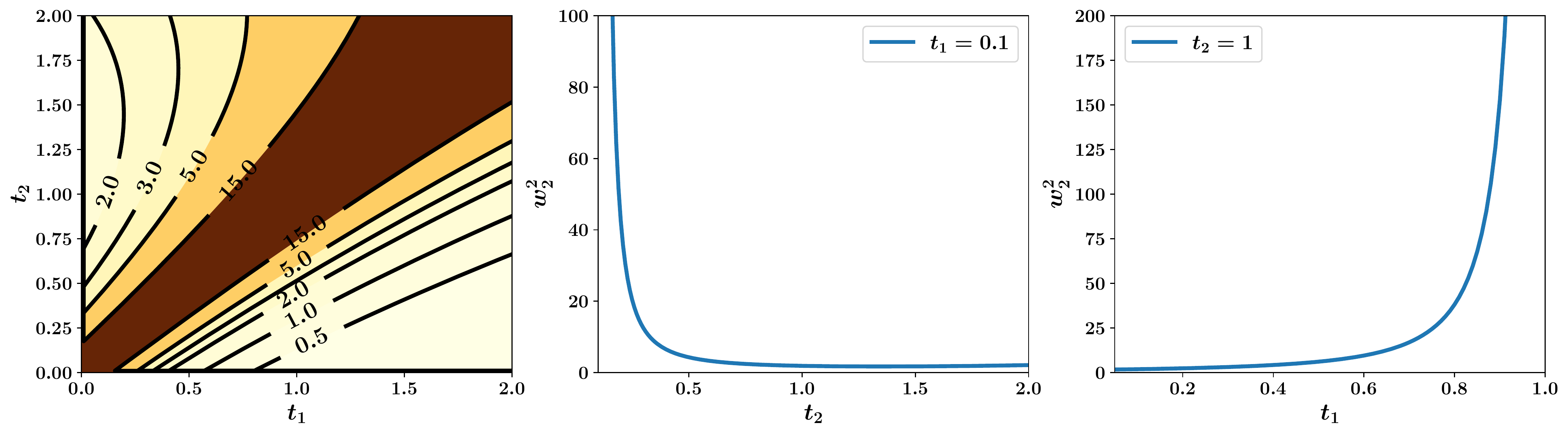}
    \caption{Sensitivity of the optimal estimates of the decay parameter $a$ to the second observation $z_2$. Left: contour plot of $w_2^2(t)$ at different values of $t_1$ and $t_2$. Middle: $y_1^2(t)$ values at $t_1=0.1$ and different values of $t_2$. Right: $w_1^2(t)$ values at $t_2=1.0$ and different values of $t_1$. Note that minimum values of $w_2^2$ are obtained at $t_1=0.1$ and $t_2=1.0$.}
    \label{fig:w2}
\end{figure}

Plots of sensitivity of the optimal control estimates with respect to observations (i.e., $y_1^2(t)$, $w_1^2(t)$, $y_2^2(t)$, and $w_2^2(t)$) are given in \cref{fig:y1} through \cref{fig:w2}. \Cref{fig:y1} shows the contour plot of $y_1^2(t)$ as well as the cross sections of this contour at the observation times (i.e., $t_1 = 0.1$ and $t_2 = 1.0$). 
Similar plots for $w_1^2(t)$ are covered in \cref{fig:w1}. It can be seen that maximum values of the sensitivity of the optimal estimate of $x_0$ with respect to the first observation occur when $t_1$ and $t_2$ are close to each other. For example, by selecting $t_1=0.1$, placing the second observation at earlier times (e.g., $t_2<0.5$) leads to very high sensitivity of the optimal estimate of $x_0$ and $a$ with respect to the first observation location. Thus, larger values of $t_2$ are to be sought to minimize $y_1^2$. 

Likewise, \cref{fig:y2} describes the properties of $y_2^2(t)$ and \cref{fig:w2} displays properties of $w_2^2(t)$. As expected from the theory, the sensitivity of the optimal control estimates with respect to the observations exhibit minimum values at the two times where the observation sites were chosen. For example, by selecting $t_1=0.1$, $y_2^2$ and $w_2^2$ approach their minimum values at $t_2=1$. In a similar way, by selecting $t_2=1$, $y_2^2$ and $w_2^2$ are minimized by placing the first observation at very small values of $t_1$.

\subsection{Nonlinear dynamics (scalar)}
\begin{figure}[ht!]
    \centering
    \includegraphics[width=0.6\linewidth]{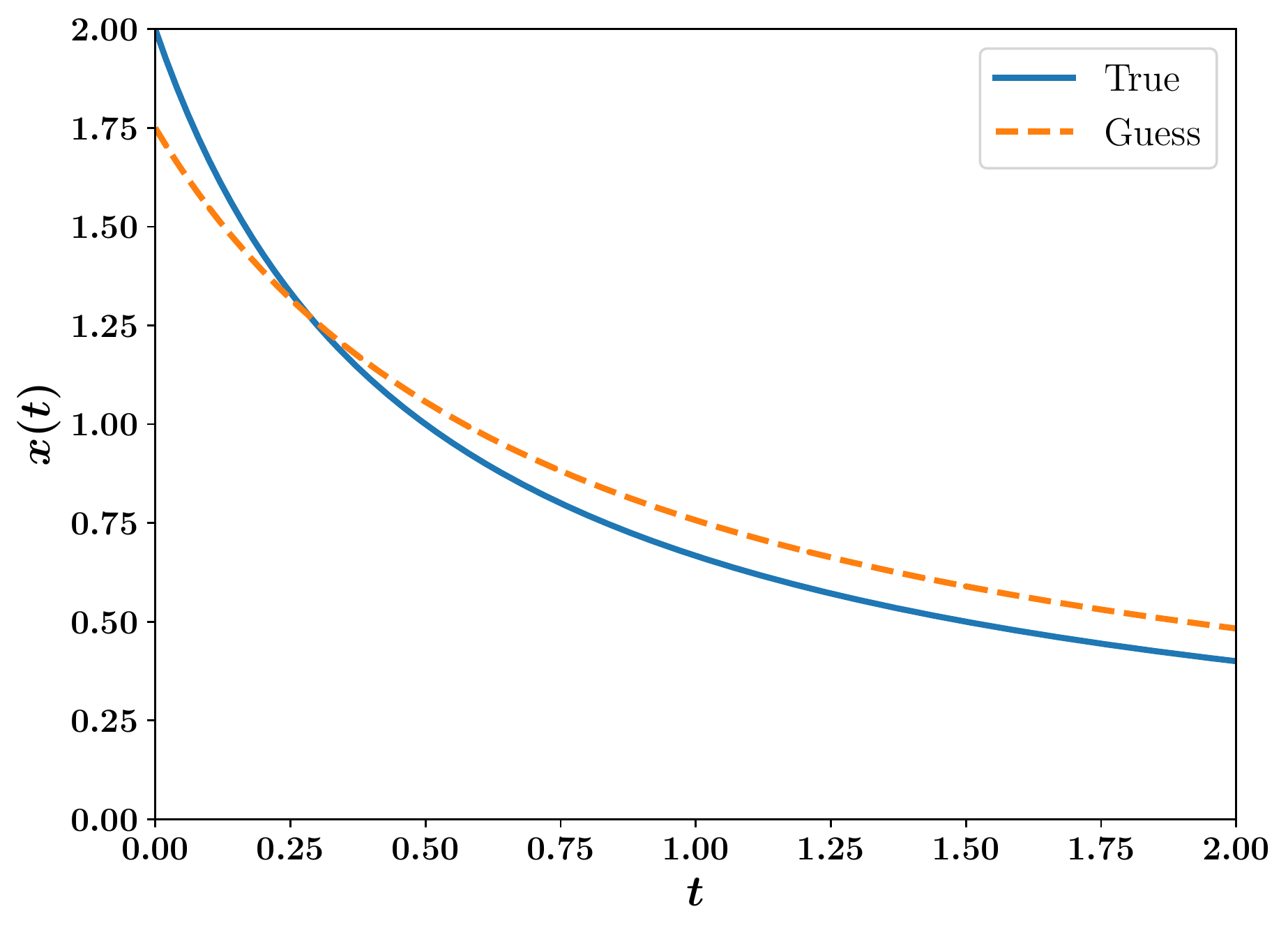}
    \caption{Solution of the nonlinear decay model using the true and guess controls.}
    \label{fig:xtruenon}
\end{figure}

We extend our analysis and consider a system governed by nonlinear dynamics as follows:
\begin{equation}
    \dot{x} = a x^2,
\end{equation}
where $x \in \mathbb{R}$ is the state variable and $a \in \mathbb{R}$ is the model parameter. It can be further shown that the solution of this system is given by $x(t) = \dfrac{x_0}{1-a x_0 t}$. The true and guessed controls are defined as $\bar{c}=[\bar{x}_0, \bar{a}]\tran = [2,-1]\tran$ and $c=[x_0, a]\tran = [1.75,-0.75]\tran$, respectively. The dynamical evolution of the nonlinear decay system in response to the true and guess controls is shown in \cref{fig:xtruenon}.

\begin{figure}[ht!]
    \centering
    \includegraphics[width=0.95\linewidth]{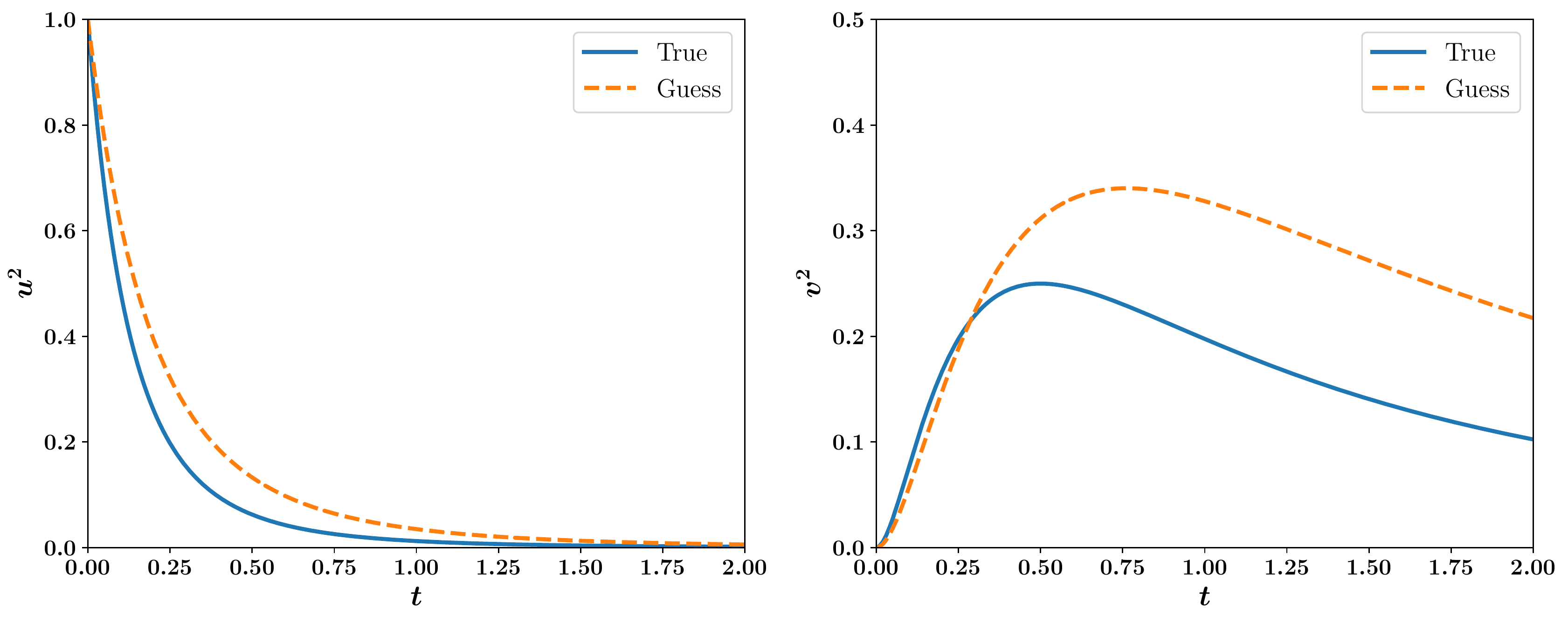}
    \caption{Model sensitivity with respect to (a) initial condition $x_0$ (left) and (b) the decay parameter $a$ (right) for $\dot{x} = a x^2$.}
    \label{fig:sensitivitynon}
\end{figure}

The sensitivity of the model forecast with respect to the initial condition and model parameter can be written as follows:
\begin{equation}
    u(t) = \dfrac{1}{(1-a x_0 t)^2}, \quad v(t) = \dfrac{x_0^2t}{(1-a x_0t)^2}.
\end{equation}
The evolution of these sensitivities with time is depicted in \cref{fig:sensitivitynon}. We note that similar figures can be obtained by solving \cref{eq:sens} numerically, where $D_f(x) = 2a x$ and $D_f^{a}(x) = x^2$.

\begin{figure}[ht!]
    \centering
    \includegraphics[width=0.65\linewidth]{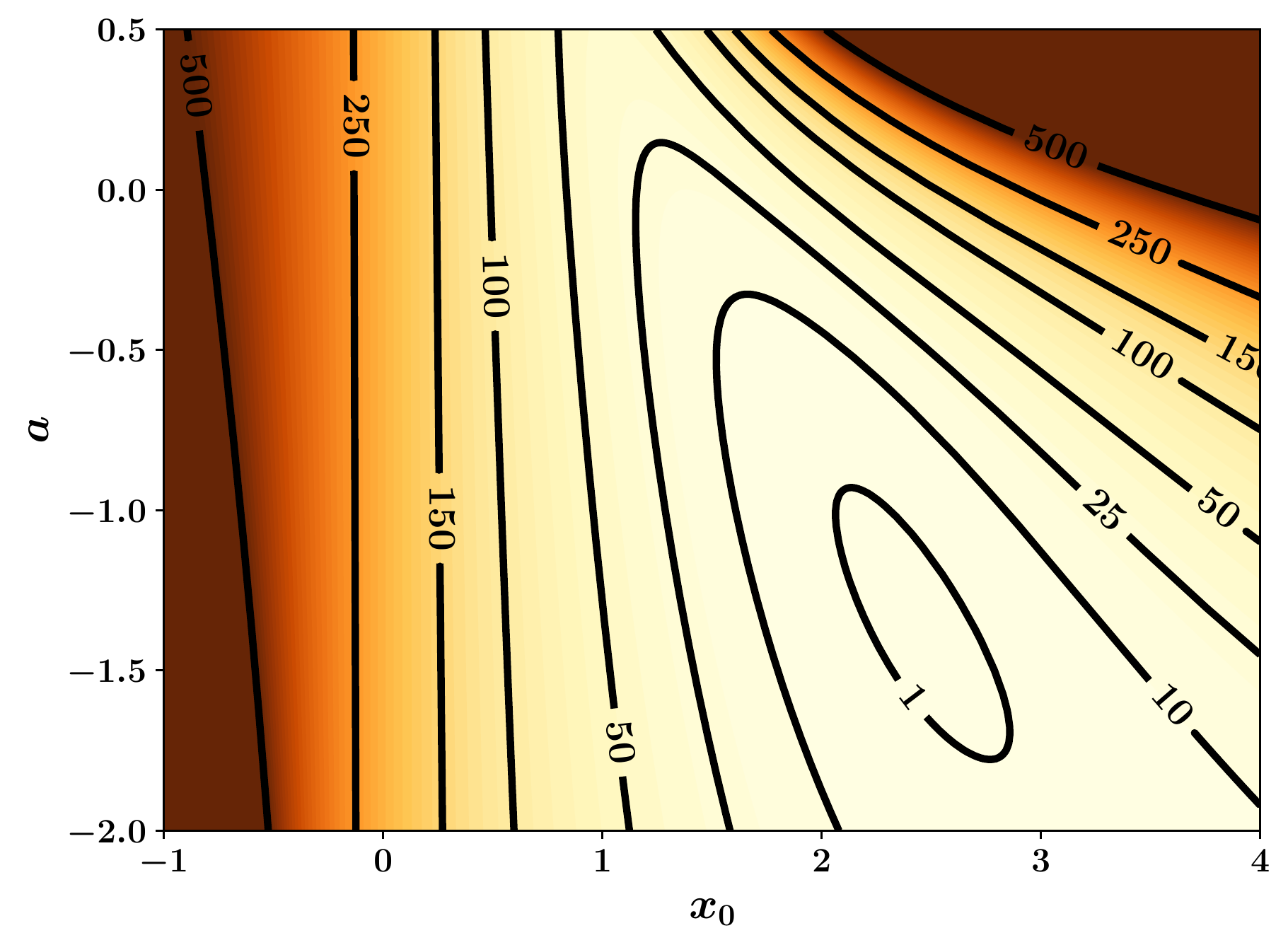}
    \caption{Cost function for nonlinear decay model: $\dot{x} = \alpha x^2$.}
    \label{fig:costnon}
\end{figure}

Although the maximum value of $u^2$ occurs at $t=0$, it corresponds to a zero value of $v^2$. Thus, we avoid placing observations at the initial time and we select $t_1=\epsilon=0.1>0$. Moreover, we select $t_2=0.5$ to maximize the value of $v^2$. We follow a twin experiment approach to generate synthetic observations by adding white Gaussian noise to the true system state (corresponding to the true controls). The resulting cost function is depicted in \cref{fig:costnon}, where Newton's method can be followed to obtain the optimal values of $x_0$ and $a$.

\begin{figure}[ht!]
    \centering
    \includegraphics[width=0.99\linewidth]{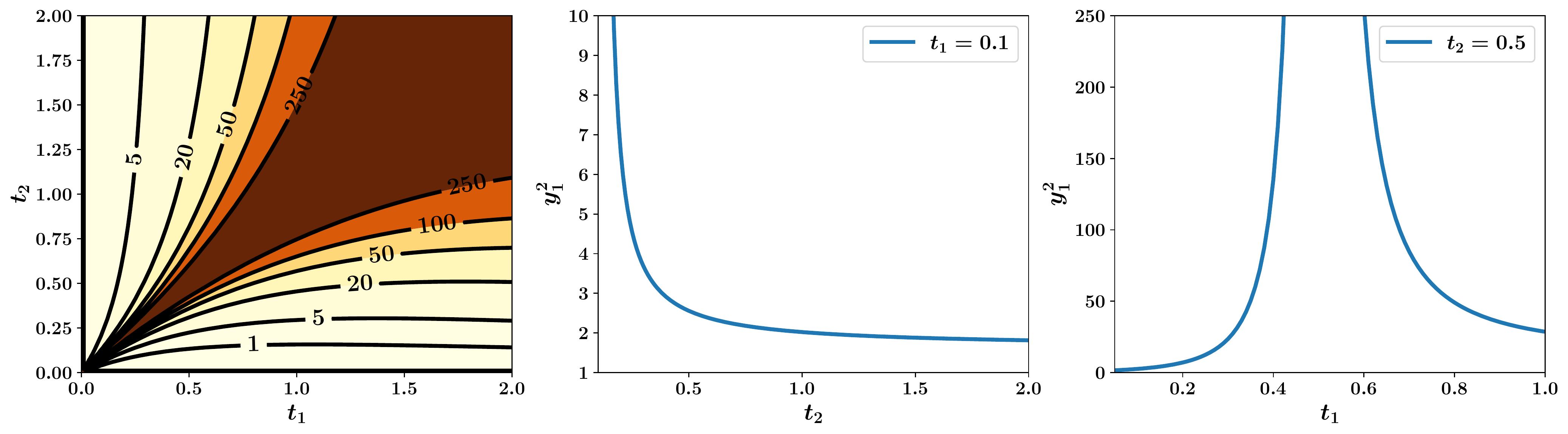}
    \caption{Sensitivity of the optimal estimates of the initial condition $x_0$ to the first observation $z_1$ for the nonlinear system $\dot{x} = \alpha x^2$. Left: contour plot of $y_1^2(t)$ at different values of $t_1$ and $t_2$. Middle: $y_1^2(t)$ values at $t_1=0.1$ and different values of $t_2$. Right: $y_1^2(t)$ values at $t_2=0.5$ and different values of $t_1$.}
    \label{fig:y1non}
\end{figure}

\begin{figure}[ht!]
    \centering
    \includegraphics[width=0.99\linewidth]{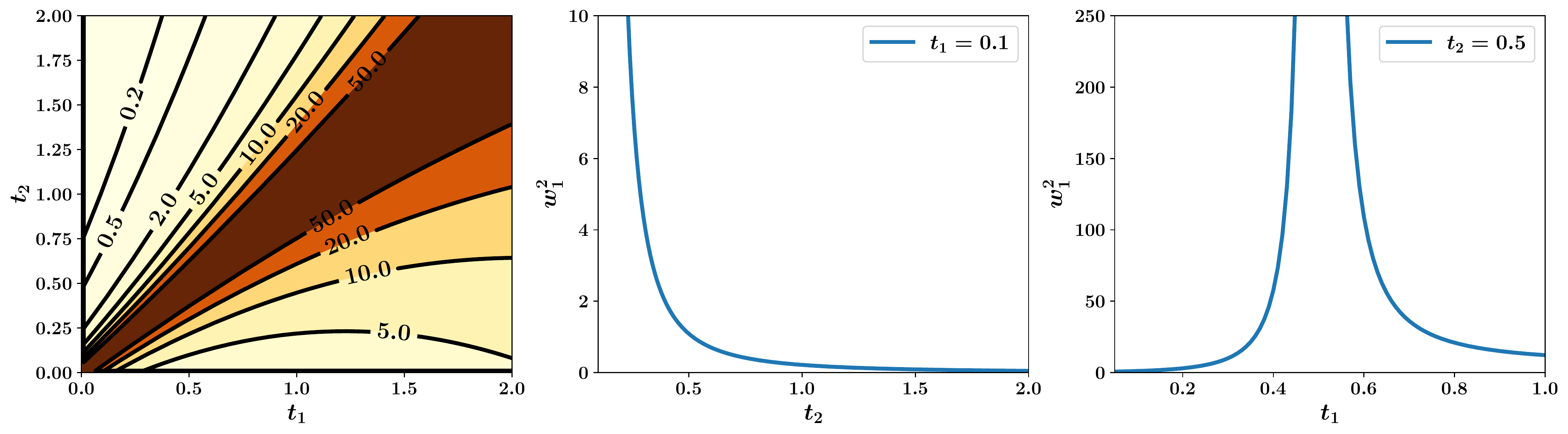}
    \caption{Sensitivity of the optimal estimates of the decay parameter $a$ to the first observation $z_1$ for the nonlinear system $\dot{x} = \alpha x^2$. Left: contour plot of $w_1^2(t)$ at different values of $t_1$ and $t_2$. Middle: $y_1^2(t)$ values at $t_1=0.1$ and different values of $t_2$. Right: $w_1^2(t)$ values at $t_2=0.5$ and different values of $t_1$.}
    \label{fig:w1non}
\end{figure}

As highlighted before, we place the observations where the squares of the forward sensitivities of the model forecast with respect to the control are maximized. This guarantees that the observability Gramian is positive definite and subsequently keeps the adjoint gradient away from zero, which in turn accelerates the search for optimal estimates of the control. The second advantage of the proposed observation placement strategy is related to the robustness of the solution of the inverse problems against small perturbations to the collected measurement data. \Cref{fig:y1non} depicts the sensitivity of the optimal estimate of initial condition with respect to the first observation, considering a total of two observations at $t_1$ and $t_2$. From the contour plots of $y_1^2$, we see that collecting the two observations at close time instants results in higher sensitivity of the estimated initial condition to the measurement itself. Similar observations are obtained from \cref{fig:w1non} for $w_1^2$, where $w_1$ denotes the sensitivity of the optimal estimate of the model parameter with respect to the first observation. This is in agreement with the fact that highly correlated measurements do not add valuable information to the inversion problem. Moreover, from \cref{eq:det}, when $t_1$ and $t_2$ coincide, the determinant of the observability Gramian vanishes, resulting in extreme values of the sensitivities $y_i$ and $w_i$ for $i=1,2$.

\begin{figure}[ht!]
    \centering
    \includegraphics[width=0.99\linewidth]{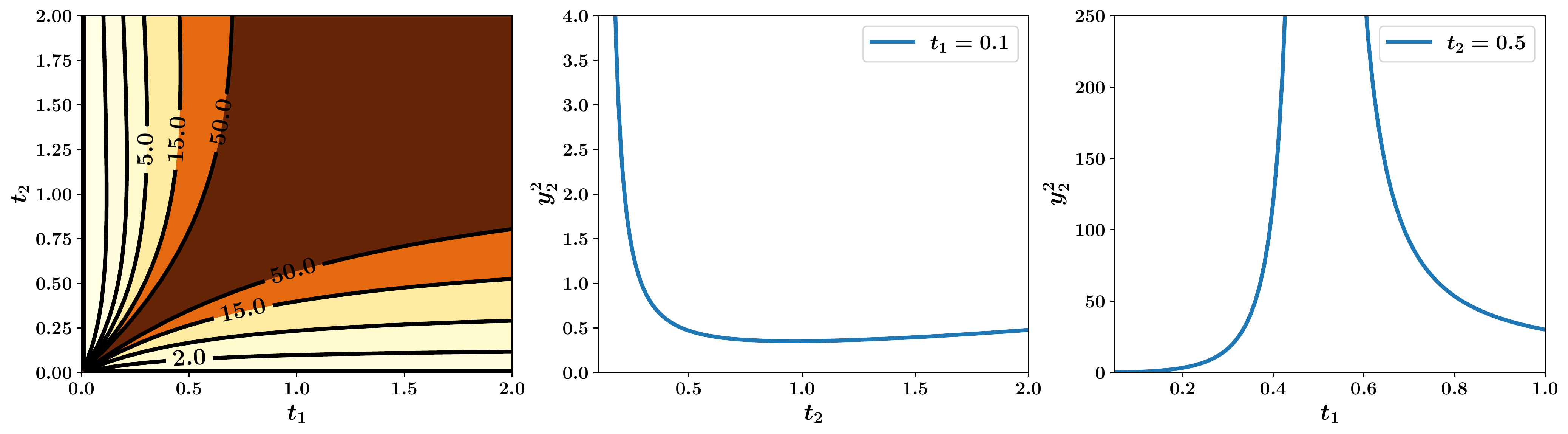}
    \caption{Sensitivity of the optimal estimates of the initial condition $x_0$ to the second observation $z_2$ for the nonlinear system $\dot{x} = \alpha x^2$. Left: contour plot of $y_2^2(t)$ at different values of $t_1$ and $t_2$. Middle: $y_2^2(t)$ values at $t_1=0.1$ and different values of $t_2$. Right: $y_1^2(t)$ values at $t_2=0.5$ and different values of $t_1$.}
    \label{fig:y2non}
\end{figure}

\begin{figure}[ht!]
    \centering
    \includegraphics[width=0.99\linewidth]{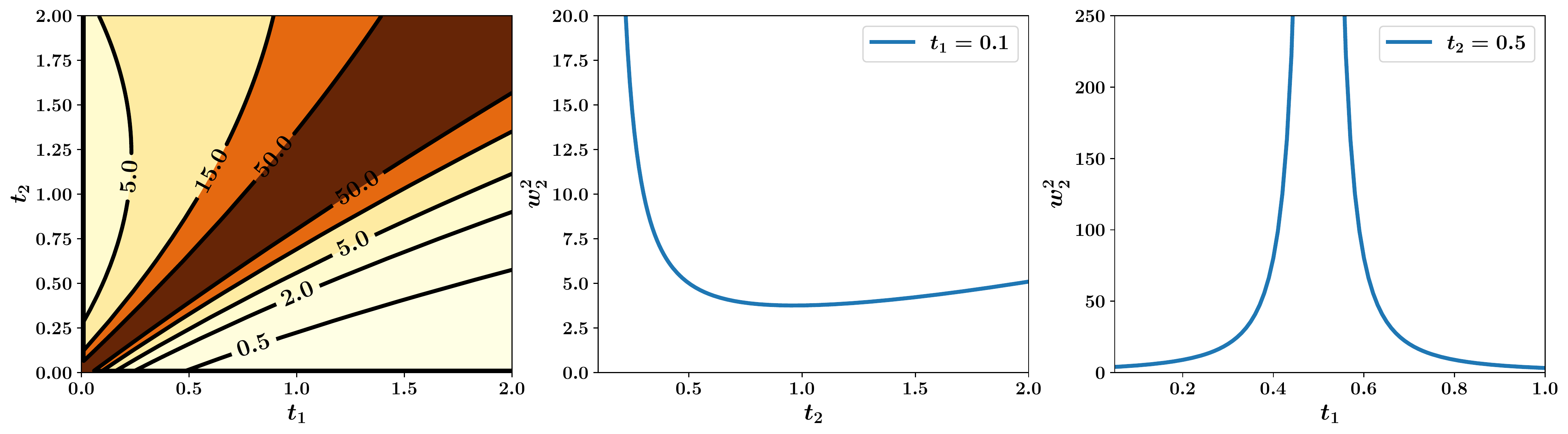}
    \caption{Sensitivity of the optimal estimates of the decay parameter $a$ to the second observation $z_2$ for the nonlinear system $\dot{x} = \alpha x^2$. Left: contour plot of $w_2^2(t)$ at different values of $t_1$ and $t_2$. Middle: $y_1^2(t)$ values at $t_1=0.1$ and different values of $t_2$. Right: $w_1^2(t)$ values at $t_2=0.5$ and different values of $t_1$.}
    \label{fig:w2non}
\end{figure}

\Cref{fig:y2non,fig:w2non} display the sensitivity of the inferred initial condition and model parameter, respectively, to the second observation value at $t_2$. Similar to $y_1^2$ and $w_1^2$ plots, we see that maximum values of $y_2^2$ and $w_2^2$ occur when $t_1$ and $t_2$ are close to each other. In addition, looking at $y_2^2$ and $w_2^2$ values at $t_1=0.1$, it is clear that $y_2^2$ and $w_2^2$ exhibit very large values at small values of $t_2$, then decrease to attain their minima around $t_2=0.5$ before slightly increasing again at larger values of $t_2$. In a similar fashion, through an investigation of $y_2^2$ and $w_2^2$ at $t_2=0.5$ with varying $t_1$, we observe that maximum values occur when $t_1=t_2=0.5$. On the other hand, minimum values of the sensitivity of the control estimate with respect to the second measurement $z_2$ at $t_2=0.5$ occur when the first measurement is collected earlier at small values of $t_1$. The numerical investigations of the sensitivity of the inverse problem solution to the observation values and locations confirm that they are minimized when the observations are placed using the proposed forward sensitivity analysis. By placing the observations at points where the squares of the forward sensitivities of the forecast are maximized, we (1) control the shape of the cost function in such a way that flat patches (with zero gradients) are avoided, which improves the convergence characteristics of solving the inverse problem, and (2) guarantee that the sensitivities of the estimated control values (solution to the inverse problem) with respect to the measurement values are minimized, which enhances the robustness of the inversion framework against possible perturbations in the collected sensor data.

\subsection{1D Burgers problem}
\begin{figure}[ht!]
    \centering
    \includegraphics[width=0.75\linewidth]{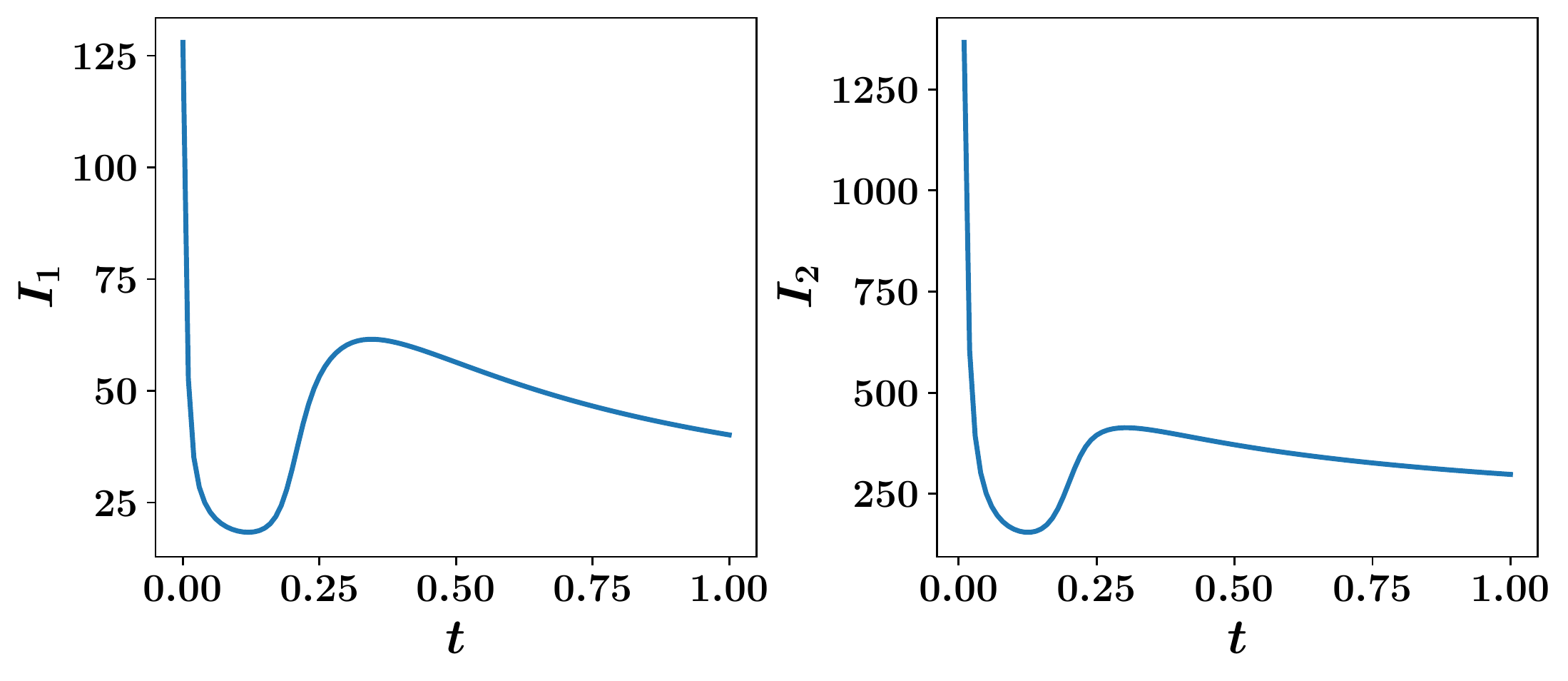}
    \caption{\textcolor{rev}{The first two invariants for the model forecast sensitivities with respect to the initial condition for the 1D viscous Burgers problem.}}
    \label{fig:brg_sens}
\end{figure}

\textcolor{rev}{Our third test case is an advection shock problem governed by the 1D viscous Burgers equation as follows:
\begin{equation}
    \dfrac{\partial u}{\partial t} + u \dfrac{\partial u}{\partial x} = \dfrac{1}{\text{Re}} \dfrac{\partial ^2 u}{\partial x^2}, \label{eq:brg}
\end{equation}
where $u(x,t)$ is the velocity field and $\text{Re}$ is the dimensionless Reynolds number, defined as the ratio of inertial effects to viscous effects. We first perform spatial discretization by defining the velocity field at $n$ discrete locations, equally spaced in the domain $L$. We apply a second-order centered finite difference scheme for the linear term and use the skew-symmetric formulation by Aref and Daripa's scheme \cite{aref1984note} for the nonlinear term as follows:
\begin{equation}
    \dfrac{\partial u_i}{\partial t} = -\dfrac{1}{3} \dfrac{(u_{i+1} + u_{i-1} + u_{i})(u_{i+1} - u_{i-1})  }{2\Delta x} + \dfrac{1}{\text{Re}} \dfrac{u_{i+1} - 2u_{i} + u_{i-1} }{\Delta x^2}, \label{eq:brg-semi}
\end{equation}
where $\Delta x$ is the grid spacing, $\Delta x = x_{i+1}-x_{i}$. By arranging the velocity field in a column vector, \cref{eq:brg-semi} is equivalent to \cref{eq:model_discrete_vector}. We define a domain of length $L=1$ and enforce zero Dirichlet boundary conditions, $u(0, t) = u(1,t) = 0$. A twin experiment is employed to generate observational data and assess the quality of the solution to inverse problems. In particular, the ground truth corresponds to $\text{Re}=500$ and the following initial condition \cite{ahmed2020long}:
\begin{equation}
    u(x, 0) = \dfrac{x}{1 + \exp{\left( \dfrac{\text{Re}}{16} \left(  4x^2 - 1 \right)  \right)} },
\end{equation}
while the background solution corresponds to a sinusoidal wave as $u(x,0) = \sin(2\pi x/L)$. Because it is unfeasible to track individual components of the forecast sensitivity matrices, we focus on the following two quantities:
\begin{equation}
    I_1 = \text{trace}\big(\mathbf{u}\tran \mathbf{u}\big), \qquad I_2 = \dfrac{1}{2} \Bigg( \bigg(\text{trace}\big(\mathbf{u}\tran \mathbf{u}\big)\bigg)^2 - \text{trace}\bigg(\big(\mathbf{u}\tran \mathbf{u}\big)^2\bigg) \Bigg).
\end{equation}
\Cref{fig:brg_sens} shows the time variation of $I_1$ and $I_2$ for the 1D Burgers problem. Both quantities attain their largest values near the initial time along with an additional \emph{bump} around $t=0.35$.
}

\begin{figure}[ht!]
    \centering
    \includegraphics[width=\linewidth]{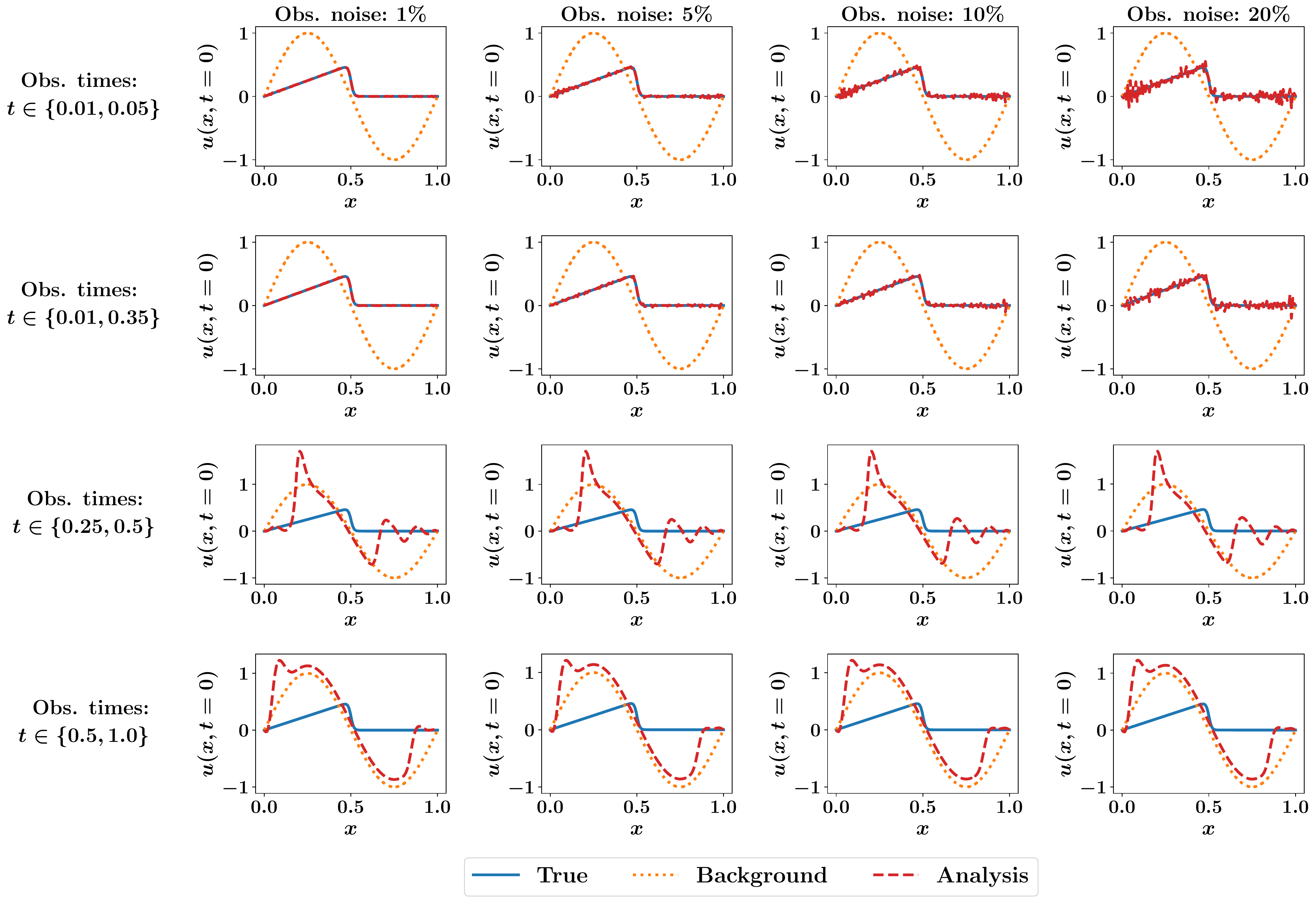}
    \caption{\textcolor{rev2}{A comparison between inverse problem solutions (for initial conditions) at different observation placement times and varying levels of noise for the 1D viscous Burgers problem. Background refers to the guessed (inaccurate) initial conditions while the analysis refers to the solution of the inverse problem after observations have been assimilated.}}
    \label{fig:brg_ic}
\end{figure}

\begin{figure}[ht!]
    \centering
    \includegraphics[width=\linewidth]{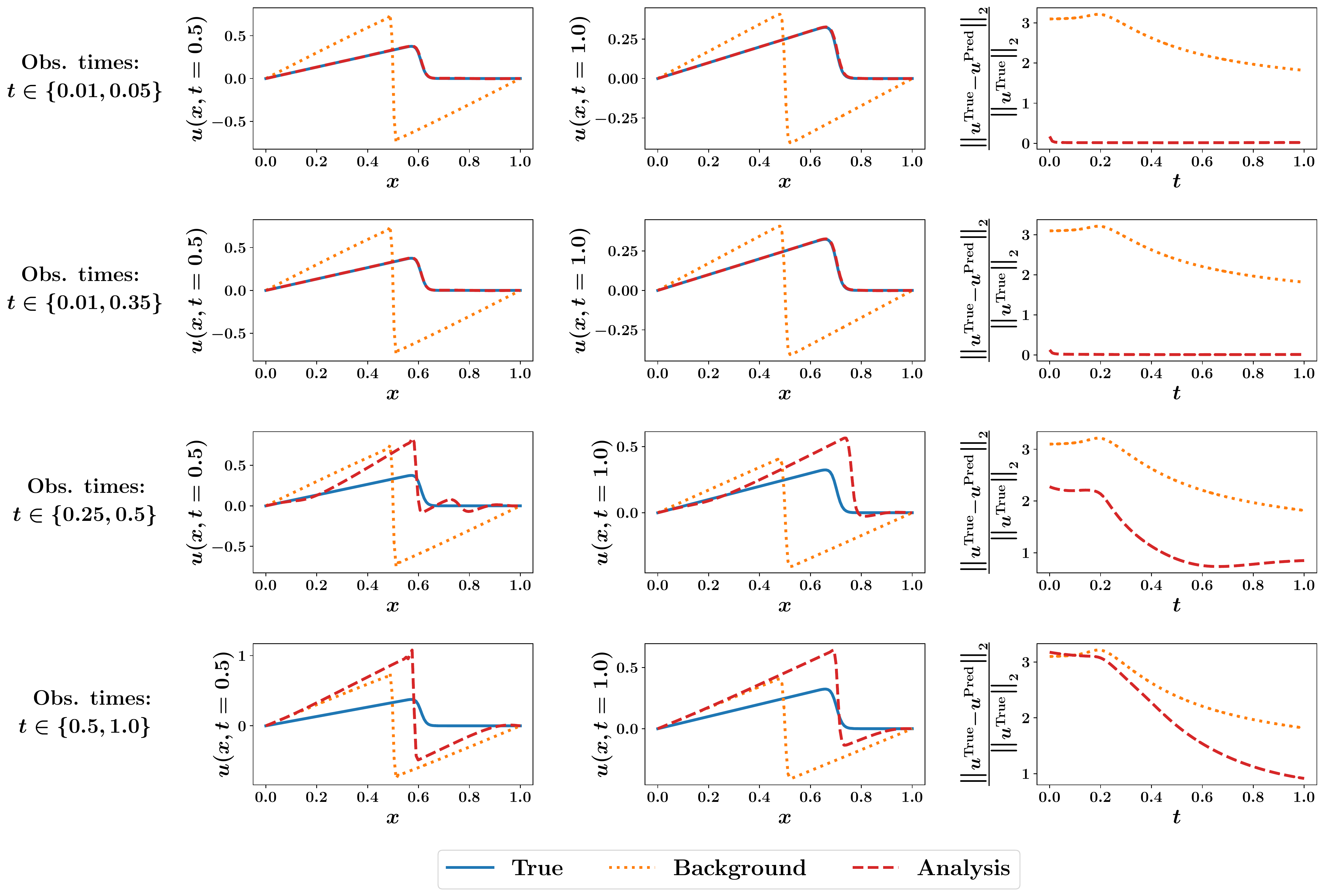}
    \caption{\textcolor{rev2}{The predicted velocity field at $t=0.5$ (left) and $t=1.0$ (middle) for the 1D viscous Burgers problem with different observation placement times with $10\%$ measurement noise along with the relative $\ell_2$ error of the predicted velocity field at different times (right). Background solution refers to the forecast using the guessed (inaccurate) initial conditions while the analysis refers to the forecast after observations have been assimilated.}}
    \label{fig:brg_pred}
\end{figure}

\textcolor{rev}{We empirically explore the advantages of the FSM-based observation placement strategy by considering different scenarios. The first two rows in \cref{fig:brg_ic} correspond to placing the observations by tracing the model forecast sensitivities. In particular, our first experiment corresponds to placing the observations near the initial time at $t=0.01$ and $t=0.05$. The second experiment corresponds to placing one observation at $t=0.01$ and another at $t=0.35$ near the second peak in \cref{fig:brg_sens}. We test the quality of the inverse problem solution with different levels of measurement noise. \Cref{fig:brg_ic} demonstrates that the FSM placement strategy enables the accurate inference of the unknown initial condition even with large levels of noise. \textcolor{rev2}{In this figure and the remaining discussion, we use standard terminology from data assimilation studies, where ``background'' refers to the estimated guess of unknown initial conditions and/or parameters. This estimate could be based on previous model predictions, historical data, or an intelligent guess. Thus, the background solution is the model forecast based on these supposedly inaccurate estimates. On the other hand, ``analysis'' refers to the solution of the inverse problem (e.g., for unknown initial conditions and/or parameters) after assimilating observation data.} We highlight that we apply a truncated singular value decomposition as a regularization technique for solving the inverse problem. Our third experiment corresponds to arbitrarily placing the observations at $t=0.25$ and $t=0.50$ while we place them at $t=0.50$ and $t=1.0$ in our fourth experiment. We can see that both cases result in an inaccurate solution of the inverse problem. It is worth noting that this inaccuracy is in part due to the off-target background information (the sinusoidal wave). However, this exaggeration is intentional to show the benefits of the proposed placement strategy even in challenging situations. Finally, the predicted velocity fields are depicted in \cref{fig:brg_pred}, where placing the observations closer to the maximum values of the forecast sensitivities yields improved results.}

\subsection{2D advection diffusion problem}
\textcolor{rev}{Our last test case is the 2D advection diffusion equation as follows:
\begin{equation}
    \dfrac{\partial u}{\partial t} + c_x \dfrac{\partial u}{\partial x} + c_y \dfrac{\partial u}{\partial y} = \nu \Big(\dfrac{\partial^2 u}{\partial x^2}+\dfrac{\partial^2 u}{\partial y^2}\Big). \label{eq:ad}
\end{equation}
\Cref{eq:ad} can be used to describe the evolution of the concentration of a substance under convection and diffusion effects, and the inversion of unknown initial conditions can be particularly important in the context of (contaminant) leak detection. Thus, variants of this test case have been used extensively in the literature of data assimilation and optimal experimental designs \cite{petra2011model,attia2022stochastic,attia2023pyoed}. We use $c_x=c_y = 0.5$ and $\nu = 0.01$ and consider an initial condition of Gaussian distribution centered at $(x_0,y_0)$ as $u(x,y,0)=e^{-\frac{1}{\nu}\big((x-x_0)^2 + (y-y_0)^2\big)}$. The true field corresponds to $(x_0,y_0)=(0.25,0.25)$ while the background (guess) solution corresponds to $(x_0,y_0)=(0.5,0.5)$.}

\begin{figure}[ht!]
    \centering
    \includegraphics[width=0.75\linewidth]{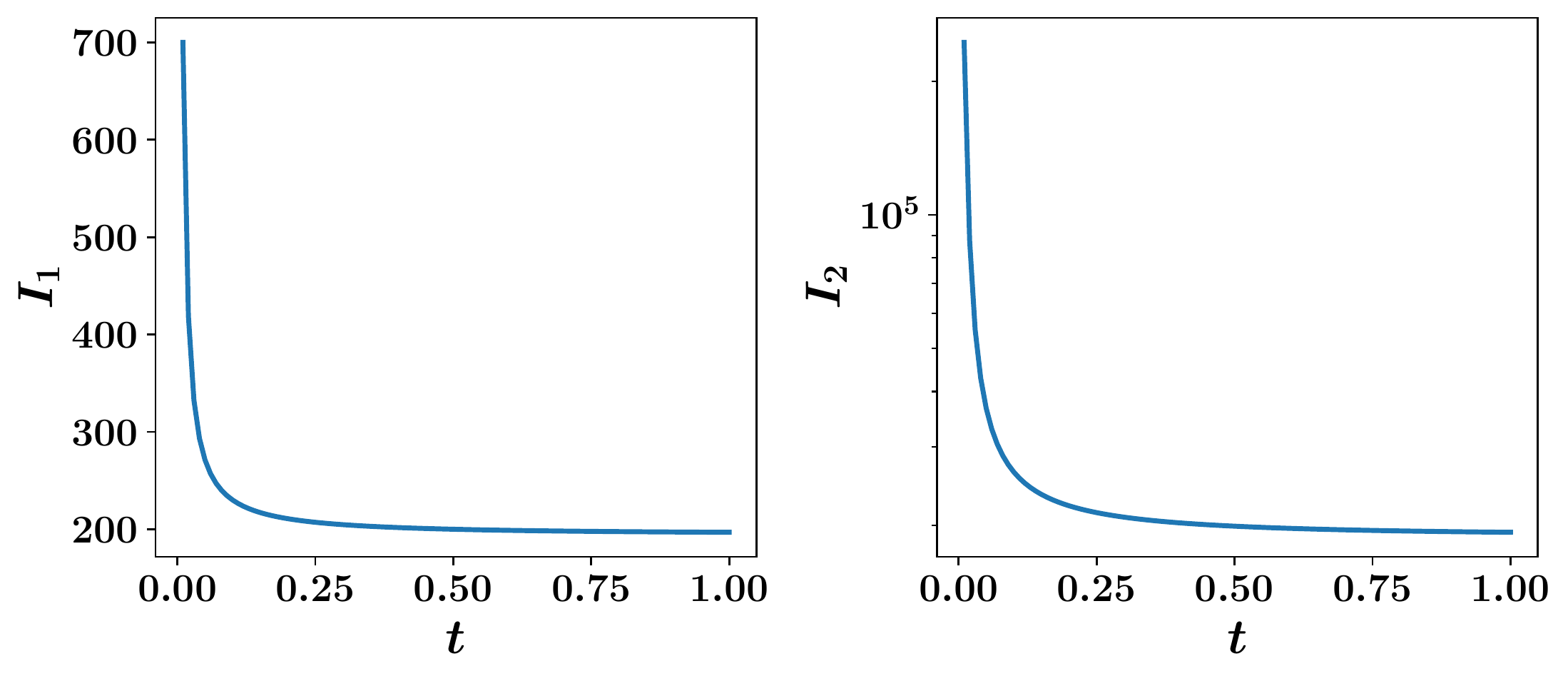}
    \caption{\textcolor{rev}{The first two invariants for the model forecast sensitivities with respect to initial conditions for the 2D advection diffusion problem.}}
    \label{fig:ad_sens}
\end{figure}

\begin{figure}[ht!]
    \centering
    \includegraphics[width=\linewidth]{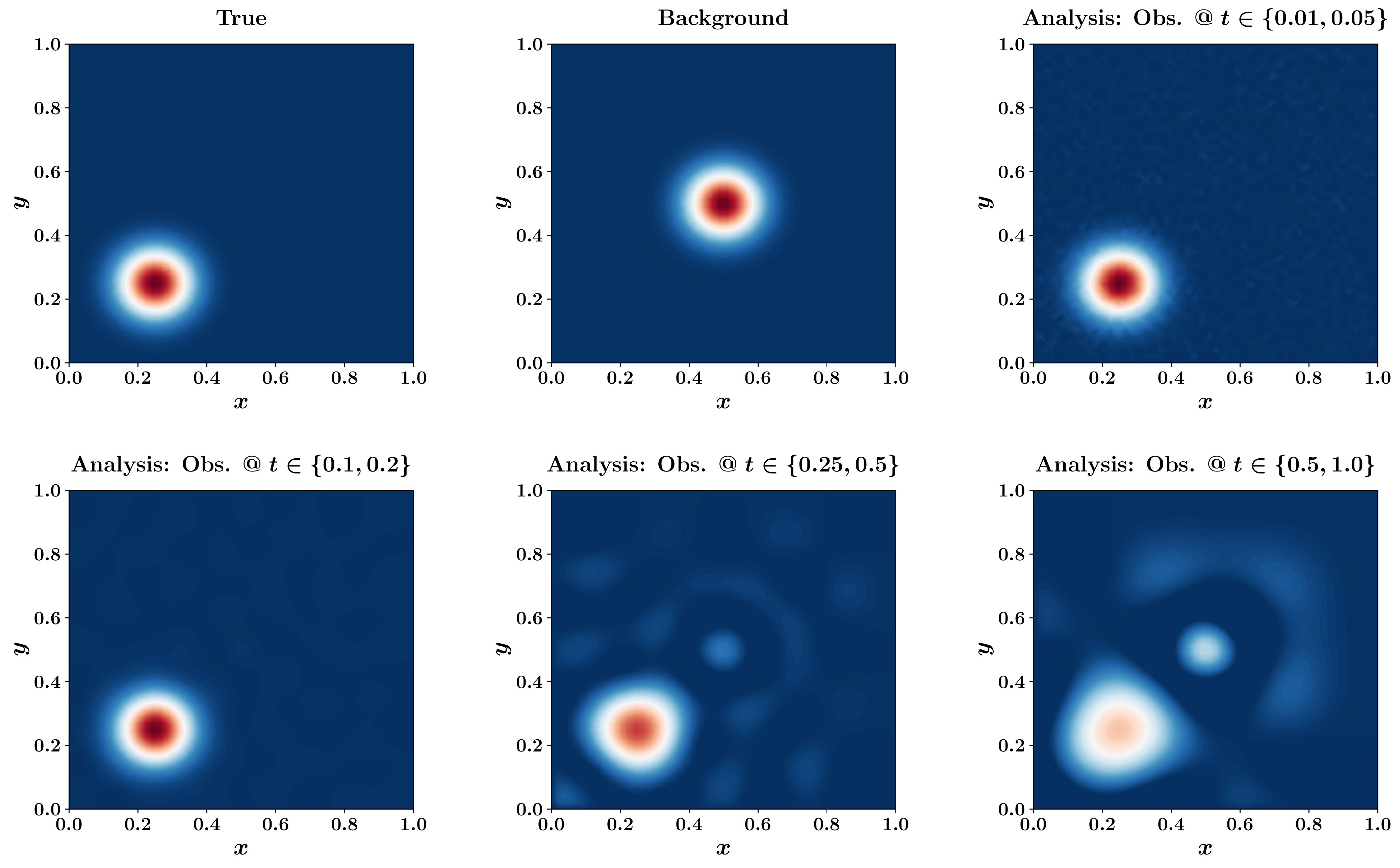}
    \caption{\textcolor{rev2}{A comparison between inverse problem solutions (for initial conditions) at different observation placement times with $10\%$ measurement noise for the 2D advection diffusion problem. Background refers to the guessed (inaccurate) initial conditions while the analysis refers to the solution of the inverse problem after observations have been assimilated.}}
    \label{fig:ad_ic}
\end{figure}

\begin{figure}[ht!]
    \centering
    \includegraphics[width=\linewidth]{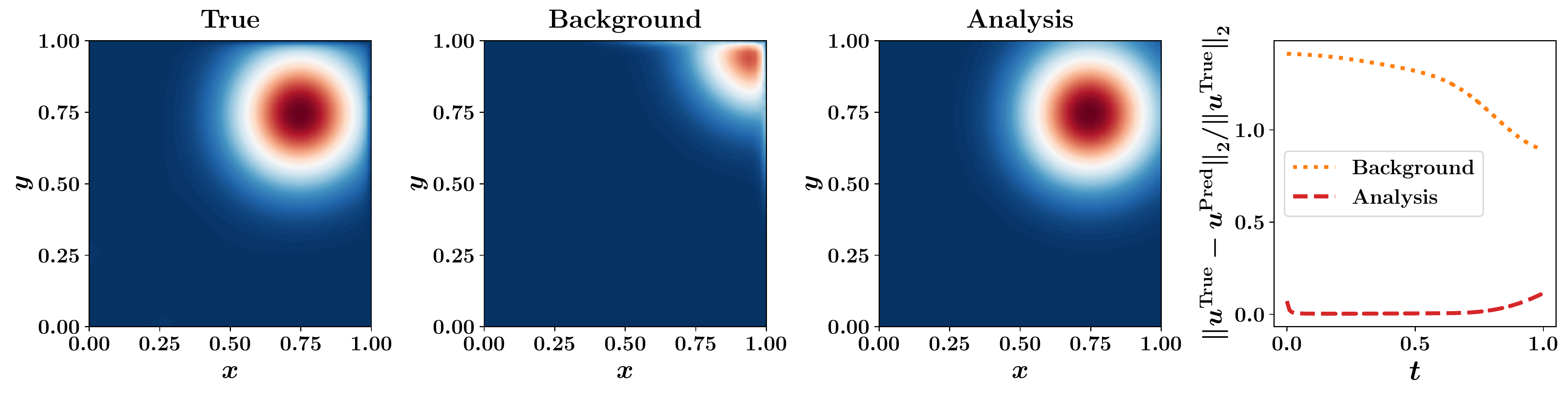}
    \caption{\textcolor{rev2}{The predicted solution at $t=1$ for the 2D advection diffusion problem following the proposed observation placement strategy (with $10\%$ measurement noise) along with the relative $\ell_2$ error of the predicted velocity field at different times (right). Background solution refers to the forecast using the guessed (inaccurate) initial conditions while the analysis refers to the forecast after observations have been assimilated.}}
    \label{fig:ad_pred}
\end{figure}

\textcolor{rev}{\Cref{fig:ad_sens} illustrates that the model forecast sensitivity with respect to the initial conditions exhibits maximum values near the initial time. In \cref{fig:ad_ic}, we repeat the solution of the inverse problems whilst varying the times at which we place our observations. While it might not be possible in practice to place all observations near the initial time, \cref{fig:ad_sens} gives theoretical guidelines to prioritize observations at specific time instants over others. For instance, it is evident from \cref{fig:ad_ic} that placing observations at $t=0.5$ and $t=1.0$ (corresponding to minimum values of forecast sensitivities) produce erroneous initial conditions. On the other hand, collecting data at $t=0.1$ and $t=0.2$ (where the forward sensitivities are relatively larger) results in a better inference of the unknown initial condition. Finally, the predicted solution at $t=1$ is shown in \cref{fig:ad_pred} where observations are placed at $t=0.01$ and $t=0.05$.}

\section{Conclusions} \label{sec:conc}
This paper systematically introduces a unified notion in the analysis and application of FSM to various questions related to forecast error correction using dynamic data assimilation \cite{lakshmivarahan2017forecast}. \textcolor{rev}{We propose a strategy based on maximizing the square of the forward sensitivities, which guarantees efficient recovery of control variables and minimizes the sensitivity of optimal control estimates with respect to observations.} First and foremost, the contribution of this paper can be better understood when compared with our ongoing efforts dedicated to the development of the FSM framework.  
In Lakshmivarahan and Lewis \cite{lakshmivarahan2010forward}, we first proved the equivalence between FSM and the 4-D VAR method. In Lakshmivarahan \emph{et al.} \cite{lakshmivarahan2020controlling,lakshmivarahan2022observability} by relating the observability Gramian to the forward sensitivities, we derived an intrinsic expression for the adjoint sensitivity in terms of the observability Gramian. By exploiting the structure of the observability Gramian, we then derived a strategy to place observations where the squares of the forward sensitivities attain maximum values. 

\textcolor{rev}{In this paper, we provide a detailed analysis of this observation placement strategy and show that it provides two crucial advantages in the context of inverse problems. First, we demonstrate that the FSM observation placement strategy avoids the occurrence of flat patches in the cost functional and improves the effectiveness of the minimization algorithm. Second, we also prove that it minimizes the sensitivity of the estimates (of the inverted quantities) with respect to the observations.} The smallness of these quantitative measurements of the sensitivity of the estimates to observations supports  determining observation sites using forward sensitivities. An advantage of the FSM analysis is that it does not require the adjoint code. However, computing the forward sensitivities is computationally demanding, but once done, we can solve two related problems. First, we can place the observations where the forward sensitivity attains a maximum value. Second, we can compute the optimal estimates of the control using FSM.

\section*{Appendix A -- Proof of the second advantage of the FSM observation placement strategy: Scalar case} \label{sec:appA} \appendix
In this Appendix, we derive an intrinsic relation that exists between the sensitivities of the estimates of the control with respect to the observations, namely $y_j$ and $w_j$ for $1\le j\le 2$, and the forward sensitivities $u(k)$ and $v(k)$.
\begin{itemize}[leftmargin=10pt]
    \item Step 1: Derivatives of $x(t)$, $u(t)$, $v(t)$, and $D_h(x(t))$ with respect to the observation $z_j$ for $1\le j\le 2$ are given by
    \begin{equation}
        \begin{aligned}
            &\dfrac{\partial x(t)}{\partial z_j} = u(t) y_j + v(t) w_j, \qquad 
             \dfrac{\partial u(t)}{\partial z_j} = \dfrac{\partial u(t)}{\partial x_0} y_j + \dfrac{\partial u(t)}{\partial \alpha} w_j, \qquad
             \dfrac{\partial v(t)}{\partial z_j} = \dfrac{\partial v(t)}{\partial x_0} y_j + \dfrac{\partial v(t)}{\partial \alpha} w_j, \\
            &\dfrac{\partial D_h(x(t))}{\partial z_j} = D_h^{(2)} (x(t)) \dfrac{\partial x(t)}{\partial z_j}, \qquad  \text{where } D_h^{(2)}(x) = \dfrac{\partial^2 h}{\partial x^2} \neq \left(\dfrac{\partial h}{\partial x}\right)^2.
        \end{aligned} \label{eq:ddz}
    \end{equation}

    \item Step 2: From the definition of $e(t)$ in \cref{eq:ek}, we get
    \begin{equation}
        \dfrac{\partial e(t)}{\partial z_j} = \delta_{tj} - D_h(x(t)) \dfrac{\partial x(t)}{\partial z_j},
    \end{equation}
    where $\delta_{tj}$ is the Kronecker delta (i.e., $\delta_{tj} = 1$ if $t=j$ and $\delta_{tj} = 0$ otherwise).
    
    \item Step 3: Differentiating $H_{1j}$ in \cref{eq:H1j} with respect to $z_p$, we get the following:
    \begin{equation}
        \left(\delta_{jp} - D_h(x(j)) \dfrac{\partial x(j)}{\partial z_p}\right) D_h(x(j)) u(j) + e(j) \dfrac{\partial D_h(x(j))}{\partial z_p} u(j) + e(j) D_h(x(j)) \dfrac{\partial u(j)}{\partial z_p} = 0 .\label{eq:H1jz}
    \end{equation}

    \item Step 4: Differentiating $H_{2j}$ in \cref{eq:H2j} with respect to $z_p$, we get the following:
    \begin{equation}
        \left(\delta_{jp} - D_h(x(j))\dfrac{\partial x(j)}{\partial z_p}\right) D_h(x(j)) v(j) + e(j) \dfrac{\partial D_h(x(j))}{\partial z_p} v(j) + e(j) D_h(x(j)) \dfrac{\partial v(j)}{\partial z_p} = 0. \label{eq:H2jz}
    \end{equation}
    
    \item Step 5: When $x(t)$ is optimal, a little reflection reveals that $e(k)$ along the optimal path is small, and henceforth it is set to zero in \cref{eq:H1jz,eq:H2jz}. Thus, the approximate expressions for $\dfrac{\partial H_{ij}}{\partial z_p}$ can be given as listed in \cref{table:dHij}.
    
    \begin{table*}[htbp!]
    \caption{Expressions for $\dfrac{\partial H_{ij}}{\partial z_p}$ for $i,j,p \in \{1,2\}$.} \vspace{5pt}
    \centering
    \begin{tabular}{P{0.05\textwidth} P{0.08\textwidth} P{0.35\textwidth} P{0.35\textwidth}   }  
    \hline
     & & $j=1$ & $j=2$ \smallskip \\
    \hline \smallskip\\
    \multirow{4}{*}{$i=1$} & $p=1$ & $\left(1 - D_h(1)\dfrac{\partial x(1)}{\partial z_1}\right) D_h(1) u(1)$ & $ - D_h^2(2) \dfrac{\partial x(2)}{\partial z_1} u(2)$  \medskip \vspace{10pt} \\
     & $p=2$ & $ - D_h^2(1) \dfrac{\partial x(1)}{\partial z_2} u(1)$ & $\left(1 - D_h(2)\dfrac{\partial x(2)}{\partial z_2}\right) D_h(2) u(2)$  \vspace{5pt} \\ 
    \hline \smallskip\\
    \multirow{4}{*}{$i=2$} & $p=1$ & $\left(1 - D_h(1)\dfrac{\partial x(1)}{\partial z_1}\right) D_h(1) v(1)$ & $ - D_h^2(2) \dfrac{\partial x(2)}{\partial z_1} v(2)$  \medskip \vspace{10pt} \\
     & $p=2$ & $ - D_h^2(1) \dfrac{\partial x(1)}{\partial z_2} v(1)$ & $\left(1 - D_h(2)\dfrac{\partial x(2)}{\partial z_2}\right) D_h(2) v(2)$  \vspace{5pt} \\ 
    \hline
    \end{tabular}
    \label{table:dHij}
    \end{table*}
    
    \item Step 6: The optimality conditions (i.e., $\dfrac{\partial H_i}{\partial z_p} = \dfrac{\partial H_{i1}}{\partial z_p} +  \dfrac{\partial H_{i2}}{\partial z_p}$) reduce to the following:
    \begin{align}
        D_h(1) u(1) &= D_h^2(1) \dfrac{\partial x(1)}{\partial z_1} u(1) + D_h^2(2) \dfrac{\partial x(2)}{\partial z_1} u(2),  \label{eq:op1}\\
        D_h(2) u(2) &= D_h^2(1) \dfrac{\partial x(1)}{\partial z_2} u(1) + D_h^2(2) \dfrac{\partial x(2)}{\partial z_2} u(2), \label{eq:op2} \\
        D_h(1) v(1) &= D_h^2(1) \dfrac{\partial x(1)}{\partial z_1} v(1) + D_h^2(2) \dfrac{\partial x(2)}{\partial z_1} v(2), \label{eq:op3} \\
        D_h(2) v(2) &= D_h^2(1) \dfrac{\partial x(1)}{\partial z_2} v(1) + D_h^2(2) \dfrac{\partial x(2)}{\partial z_2} v(2). \label{eq:op4}
    \end{align}
    
    \item Step 7: Substituting for $\dfrac{\partial x(t)}{\partial z_j}$ using \cref{eq:ddz} in \cref{eq:op1,eq:op3}, we get
    \begin{align}
        D_h(1) u(1) = D_h^2(1) \left( u^2(1) + u(1)v(1) \right) y(1) + D_h^2(2) \left( u^2(2) + u(2)v(2) \right) w(1), \label{eq:op11} \\ 
        D_h(1) v(1) = D_h^2(1) \left( u(1)v(1) + v^2(1) \right) y(1) + D_h^2(2) \left( u(2) v(2) + v^2(2) \right) w(1). \label{eq:op33}
    \end{align}
    
\end{itemize}

Referring to \cref{eq:F,eq:G}, we define $F(i) = [u(i), v(i)]\tran \in \mathbb{R}^{2\times1}$ and $G(i) = F\tran(i) D_h^2(i) F(i) \in \mathbb{R}^{2\times2}$. Thus, \cref{eq:op11} and \cref{eq:op33} can be rewritten as 
\begin{equation}
    D_h(1) \begin{bmatrix} u(1) \\ v(1) \end{bmatrix} = G \begin{bmatrix} y_1 \\ w_1 \end{bmatrix}, 
\end{equation}
where $G = G(1) + G(2)$. Therefore, $y_1$ and $w_1$ can be obtained as
\begin{equation}
    \begin{bmatrix} y_1 \\ w_1 \end{bmatrix} = G^{-1} D_h(1) \begin{bmatrix} u(1) \\ v(1) \end{bmatrix}. \label{eq:y1w1}
\end{equation}
Similarly, from \cref{eq:op2,eq:op4} and using the same procedure, we get 
\begin{equation}
    \begin{bmatrix} y_2 \\ w_2 \end{bmatrix} = G^{-1} D_h(2) \begin{bmatrix} u(2) \\ v(2) \end{bmatrix}. \label{eq:y2w2}
\end{equation}

\section*{Appendix B -- Linear independence of the forward sensitivities of the model forecast at different times with respect to the control} \label{sec:appB} \appendix
The discrete-time version of the model given in \cref{eq:model} can be written as:
\begin{equation}
    x_{k+1} = M(x_k,\alpha), \label{eq:model_discrete}
\end{equation}
where $M:\mathbb{R} \times \mathbb{R} \to \mathbb{R}$ defines the one-time step mapping defined by applying a temporal integration scheme. \Cref{eq:sens} for the dynamics of $u(t)$ and $v(t)$ can be rewritten as follows:
\begin{equation}
\begin{aligned}
    u_{k+1} &= D_M(k) u_k, \quad u_0 = 1,\\
    v_{k+1} &= D_M(k) v_k + D_M^{\alpha}(x_k), \quad v_0 = 0,
\end{aligned} \label{eq:sens_discrete}
\end{equation}
where
\begin{equation}
    D_M(k) = \dfrac{\partial M}{\partial x}\bigg|_{x=x_k}, \quad \text{and} \quad D_M^{\alpha}(k) = \dfrac{\partial M}{\partial \alpha}\bigg|_{x=x_k}.
\end{equation}

For simplicity of notation, let $a_k:=D_M(k)$ and $b_k:=D_M^{\alpha}(k)$. Thus, the sequence of $u_k$ can be a written as follows:
\begin{equation}
    \begin{aligned}
        u_k & = a_{k-1} a_{k-2} \dots a_1 a_0 , \qquad \text{for } \quad k \ge 1.
    \end{aligned}
\end{equation}
Similarly, the sequence of $v_k$ can be expanded as follows:
\begin{equation}
    \begin{aligned}
        v_k &= (a_{k-1} a_{k-2} \dots a_1 b_0) 
             + (a_{k-1} a_{k-2} \dots a_2 b_1)
             + (a_{k-1} a_{k-2} \dots a_3 b_2) 
             + \dots \\
             &+ (a_{k-1} a_{k-2} b_{k-3}) + (a_{k-1} b_{k-2}) + (b_{k-1}), \qquad \text{for } \quad k \ge 1. 
    \end{aligned}
\end{equation}

For the Gramian matrix $G$ in \cref{eq:gramian2} to be singular (i.e., the determinant in \cref{eq:det} equals zero), the vectors $[u_1,u_2]\tran$ and $[v_1,v_2]\tran$ ought to be linearly dependent. In what follows, we show that this situation cannot happen. For the vectors $[u_1,u_2]\tran$ and $[v_1,v_2]\tran$ to be linearly dependent, we get $u_k=C v_k$ for all $k$ where $C$ is a non-zero constant, then we can write the following:
\begin{align}
    u_1 &= C v_1 \Longleftrightarrow a_0 = C b_0, \\
    u_2 &=  C v_2 \Longleftrightarrow a_1 a_0 = Ca_1 b_0 + C b_1. \label{eq:u2Cv2}
\end{align}
Because $a_0 = C b_0$, \cref{eq:u2Cv2} can be rewritten as  $a_1 a_0 = a_1 a_0 + C b_1$, leading to $b_1 = 0$. Similarly, 
\begin{align}
    u_3 &= C v_3  \Longleftrightarrow a_2 a_1 a_0 = Ca_2 a_1 b_0 + C a_2 b_1 + C b_2. 
\end{align}
With $a_0 = C b_0$ and  $b_1 = 0$, we also get $b_2 = 0$. Following the same procedure, this leads to $b_k=D_M^{\alpha}(k)=0$ for all $k\neq 0$. In other words, this implies that the sensitivity of the model $M$ with respect to the parameter $\alpha$ equals zero along the trajectory except at $t=0$. However, this cannot be the case because a model cannot be sensitive to its parameters at $t=0$ and not sensitive at other times. Therefore, the vectors of $u_k$ and $v_k$ cannot be linearly dependent, and the proposed observation placement strategy thus does not lead to a singular observability Gramian.

\section*{Acknowledgments}
O.S. was supported by the U.S. Department of Energy (DOE), Office of Science, Advanced Scientific Computing Research (ASCR) program under Award Number DE-SC0019290. The work of S.E.A is supported by the DOE ASCR program through the Pacific Northwest National Laboratory Distinguished Fellowship in Scientific Computing (Project No. 71268). Pacific Northwest National Laboratory is operated by Battelle Memorial Institute for DOE under Contract DE-AC05-76RL01830.

\section*{Data Availability}
The data that support the findings of this study are available within the article.

\bibliographystyle{unsrt} 

\bibliography{manuscript}

\end{document}